 \documentclass[11pt]{article}
 \usepackage{stmaryrd}
 \usepackage{amsmath}
 \usepackage{pifont}
 \usepackage{amsfonts}
 \usepackage{mathrsfs}

 \usepackage{booktabs}
 \usepackage{threeparttable}
 \usepackage{mathrsfs,amsfonts,amsmath}
 \usepackage[dvips]{color}

\newcommand\email[1]{\href{mailto:#1}{ \nolinkurl{#1}}}

 \newtheorem{theorem}{Theorem}[section]
 \newtheorem{definition}[theorem]{Definition}
 \newtheorem{lemma}[theorem]{Lemma}
 \newtheorem{corollary}[theorem]{Corollary}
 \newtheorem{proposition}[theorem]{Proposition}
 \newtheorem{remark}[theorem]{Remark}
 \newtheorem{condition}[theorem]{Condition}
 \newtheorem{example}{Example}[section]

 \def\blemma{\begin{lemma}}\def\elemma{\end{lemma}}
 \def\bproposition{\begin{proposition}}\def\eproposition{\end{proposition}}
 \def\ttheorem{\begin{theorem}}\def\etheorem{\end{theorem}}
 \def\bcorollary{\begin{corollary}}\def\ecorollary{\end{corollary}}
 \def\bremark{\begin{remark}}\def\eremark{\end{remark}}
 \def\bcondition{\begin{condition}}\def\econdition{\end{condition}}

 \def\benumerate{\begin{enumerate}}\def\eenumerate{\end{enumerate}}
 \def\bitemize{\begin{itemize}}\def\eitemize{\end{itemize}}

 \def\beqlb{\begin{eqnarray}}\def\eeqlb{\end{eqnarray}}
 \def\beqnn{\begin{eqnarray*}}\def\eeqnn{\end{eqnarray*}}
 \def\ar{\!\!\!&}

 \def\mcr{\mathscr}\def\mbb{\mathbb}\def\mbf{\mathbf}

 \def\proof{\noindent{\it Proof.~~}}\def\qed{\hfill$\Box$\medskip}

\begin{document}

 \title{\bf\Large BACKWARD DOUBLY STOCHASTIC EQUATIONS WITH JUMPS AND COMPARISON THEOREMS\thanks{Supported by NSFC (No. 11401012).}
 \footnotetext[2]{\textit{MSC2010 subject classifications:} 		60H05, 60H10, 37H10}
 }

 \author{\large\scshape  Wei Xu\\[1mm]
 \textit{ Beijing Normal University}}
 \maketitle

 \renewcommand{\abstractname}{}
\begin{abstract}
\centerline{\scshape\large Abstract}
\medskip

In this work we mainly prove the existence and pathwise uniqueness of solutions to general backward doubly stochastic differential equations with jumps appearing in both forward and backward integral parts.  Several comparison theorems under some weak conditions are also given. Finally we apply comparison theorems in proving the existence of solution to some special backward doubly stochastic differential equations with drift coefficient increasing linearly.

\end{abstract}

  \section{Introduction}
 \setcounter{equation}{0}
 \medskip

 Backward stochastic differential equations (BSDEs) in the linear case
 were introduced in Kushner (1972), Bismut (1976), Bensoussan (1982) and Haussmann (1986) as adjoint processes in the
 maximum principle for stochastic control problems and the pricing of options.
 Since the  important work of Pardoux and Peng (1990), the interest in BSDEs has increased considerably in recent years.
 The significance of BSDEs is not only proved by the considerably important role they are playing in the study of partial differential equations (PDEs); see Peng (1991), Pardoux and Peng (1992) and Darling and Pardoux (1997), but also can be found in many other fields such as mathematical economics, financial mathematics, insurance and stochastic control. Here we just list  several important works in every field. Duffie and Epstein (1992a,b) used BSDEs as a powerful tool to study stochastic differential utility.
 Moreover, in the insurance market BSDEs are used in
 pricing and hedging insurance equity-linked claims and asset-liability management problems, see El Karoui et al. (1997) and Delong (2013).
 Peng (1993) studied stochastic optimal control systems, where the state variables are described by a system of ordinary-SDE and BSDEs, and derived a local form of the maximum principle.

 As further extensions of BSDEs, backward doubly stochastic differential equations (BDSDEs) contain both forward and backward stochastic integrals.
 Those equations were first introduced by Pardoux and Peng (1994) in the study of quasi-linear parabolic stochastic partial differential equations (SPDEs).
 Compared to BSDEs,  much less results about BDSDEs can be found in the literature and most of the results established are about BDSDEs driven by Brownian motions.
 For the details about applications of BDSDEs to SPDEs driven by Brownian motion, one can refer to Zhang and Zhao (2007) which studied the existence and uniqueness of solution to BDSDEs on infinite horizons, and the stationary solutions to SPDEs by virtue of the solutions to BDSDEs on infinite horizons.
 Moreover, some work about BDSDEs with jumps appearing in the system of ordinary-SDE have been published recently.
 For instance, Zhu and Shi (2012) studied BDSDEs driven by Brownian motions and Poisson process with non-Lipschitz coefficients on random time interval.
 Aman (2012), Aman and Owo (2012) and Ren et al. (2009) study a special reflected generalized BDSDEs (driven by Teugel's martingales associated with L\'{e}vy process) with means of the penalization method and the fixed-point theorem. Existence and uniqueness of the solution to the BDSDE with jumps in the forward integral are studies in Sow (2011) for the case of non-Lipschitz coefficients.
 Recently, some results about stochastic control problems of BDSDEs have been obtained by Han et al. (2010) and Bahlali and Gherbal (2010).

 This work is motivated by Xiong (2013) and He et al. (2014) which mainly studied the distribution function valued process of super-Brownian motions and super-L\'evy processes characterized as the pathwise unique solution to a SPDE.
 For any super-L\'evy process with transition semigroup $(Q_t)_{t\geq 0}$ defined by (1.4) in He et al. (2014), they proved that its distribution function valued process solves the following stochastic integral equation: for any $x\in \mathbb{R}$,
  \beqlb\label{1.1}
 X_t(x)\ar=\ar X_0(x) +\int_0^t A^*X_s(x)ds +\sqrt{c}\int_0^t\int_0^{X_{s-}(x)} W(ds,du)\cr
 \ar\ar +\int_0^t\int_0^{\infty}\int_0^{X_{s-}(x)}z\tilde{N}(ds,dz,du)-b\int_0^tX_s(x)ds,
 \eeqlb
 where $\{W(ds,du);t\geq 0, u >0\}$ is a Gaussian white noise with intensity $ds\pi(du)$,  $\{N (dt,dz,du): t\geq 0, u >0\}$ is a Poisson random measure with intensity $dt\mu(dz)du$ and $A^*$ is the dual operator of $A$ defined by: for any $f(x)\in C_0^2(\mathbb{R})$,
 \beqnn
 Af(x)=\beta f'(x)+\frac{1}{2}\sigma^2f''(x) +\int_{\mathbb{R}}[f(x+z)-f(x)-f'(x)z\mathbf{1}_{\{|z|\leq 1\}}]\nu(dz).
 \eeqnn
 Furthermore, for any fixed $T>0$, define the
 Gaussian white noise $\{W^T(dt,du): 0\le t\le T, u \in E\}$ by
 $$
 W^T([T-t, T ]\times A) = W([0, t]\times A), \qquad 0\le t\le T, A \in \mcr{B}(E);
 $$
 and the Poisson random measures $\{N^T(dt,du):0\le t\le
 T, u \in U_0\}$ by:
 $$
 N^T([T-t, T ]\times B)= N ( [0, t]\times B), \qquad 0\le t\le T, B \in \mcr{B}(U_i ).
 $$
 In the proof of the pathwise uniqueness of solutions to (\ref{1.1}), they established its connection to the following BDSDE:
 \beqnn
 X_{T-t}(L_t^r+x) \ar=\ar X_0(L_t^r+x)-b\int_t^TX_{T-s}(L_s^r+x)ds +\sqrt{c}\int_{t-}^{T-}\int_0^{X_{T-s}(L_s^r+x)} W(ds,du)\cr
 \ar\ar +\int_{t-}^{T-}\int_0^{\infty}\int_0^{X_{T-s}(L_s^r+x)}z\tilde{N}_0^T(\overleftarrow{ds},dz,du)-\sigma\int_t^T \nabla X_{T-s}(L_s^r+x)dB_s\cr
 \ar\ar -\int_t^T \nabla [X_{T-s}(L_s^r+x-z)-X_{T-s}(L_s^r+x)]\tilde{M}(ds,dz),
 \eeqnn
 where $L_t^r=L_t-L_r$ and $\{L_t:t\geq 0\}$ is a L\'evy process with generator $A^*$.

 The purpose of this work is extending the above equations into more general BDSDEs with jumps appearing not only in the forward stochastic integral part but also in the backward stochastic integral part; see (\ref{2.1}) in Section~2.  Pathwise uniqueness and existence of their solutions are proved in Section~2 and 3 respectively under Lipschitz conditions.
 In addition, several comparison theorems for BDSDEs will also be given in Section 4, since they play an important role in both theory and applications; see Shi et al. (2005).
 Effected by random terms in the backward integrals, classical methods are not applicable,
 we use another method to get comparison theorems with some reasonable and weak conditions.
 The main difficulty is to deal with the influence of forward integrals to the drift coefficient and backward integrals.
 As an applications of comparison theorems, in Section~5 we prove that solutions to a special kind of BDSDEs with drift coefficient increasing linearly exist.

 \textbf{Notation:}
 For any $n$-dimensional vector $X=(x_1,\cdots,x_n)$, $Y=(y_1,\cdots,y_n)$ and $n\times n$-matrix $A=(a_{i,j})$, let $\|X\|^2=\sum_{i=1}^n x_i^2$, $\langle X,Y\rangle=\sum_{i=1}^n x_iy_i$, $\mathbf{T}(A)=(a_{11},\cdots,a_{nn})$ and $\|A\|^2={\rm Tr}(A^{\rm T}A)=\sum_{i,j=1}^n a_{ij}^2$, where ${\rm Tr}(A)$ is the trace of $A$. For any $f\in C^2(\mathbb{R}^n)$, let
 $$Df(x)=\Big(\frac{\partial f(x)}{\partial x_1},\cdots,\frac{\partial f(x)}{\partial x_n}\Big)\quad \mbox{and}\quad D^2f(x)=\Big(\frac{\partial^2 f(x)}{\partial x_1^2},\cdots,\frac{\partial^2 f(x)}{\partial x_n^2}\Big).$$
  Throughout this paper, we make the conventions
 \beqnn
 \int_a^b\! = \int_{(a,b]}\quad,\quad \int_{a}^{\infty}\! = \int_{(a,\infty)}\quad \mbox{and}\quad \int_{a-}^{b-}\! = \int_{[a,b)}
 \eeqnn
 for any $b \geq a \geq  0$.

\section{Pathwise Uniqueness}
\setcounter{equation}{0}
 \medskip

  In this section, we mainly study the pathwise uniqueness of solutions to general backward doubly-stochastic equations.
  Suppose that $T> 0$ is a fixed constant and $(\Omega,\mcr{F}, \mbf{P})$ is a complete probability space endowed with filtration
  $\{\mcr{G}^{(1)}_t\}_{0\leq t\leq T}$ satisfying the usual hypotheses.
  Let $B(s)$ is a $n$-dimensional $(\mcr{G}^{(1)}_t)$-Brownian motion and $\{M(dt,du):0\le t\le T, u\in F\}$ a $(\mcr{G}^{(1)}_t)$-Poisson random measure with intensity $dt\nu(du)$, where $\nu(du)$ is a $\sigma$-finite Borel measure on the Polish spaces $F$.
  Let $\{\mcr{G}^{(2)}_t\}_{0\leq t\le T}$ be another filtration on $(\Omega,\mcr{F}, \mbf{P})$ satisfying the usual hypotheses and independent from $\{\mcr{G}^{(1)}_t\}_{0\leq t\le T}$. Let $\{W(ds,du)=(W_1(ds,du),\cdots,
  W_n(ds,du))^{\rm T};0\le t\le T, u \in
  E\}$ be a $n$-dimensional $(\mcr{G}^{(2)}_t)$-Gaussian white noise constructed with $n$ orthogonal white noises $W_i(ds,du)$ on $\mathbb{R}^+\times E$ with intensity $ds\pi_i(du)$ respectively. Here we denote $\pi(du)=(\pi_1(du),\cdots,\pi_n(du))^{\rm T}$.
  Suppose $\mu_0(du)$ is a $\sigma$-finite Borel measure on the
  Polish space $U_0$ and $\mu_1 (du)$ is
  a finite Borel measure on the Polish space $U_1$.  Moreover, For each $i = 0, 1$, let $\{N_i (dt, du): 0\le t\le T, u \in
  U_i\}$ be a $(\mcr{G}_t^{(2)})$-Poisson random measure with intensity
  $dt\mu_i (du)$. Obviously, all the random elements introduced
  above are independent of each other.

  Denote $\mcr{G}_t^r = \sigma(\mcr{G}^{(1)}_t\cup
  \mcr{G}^{(2)}_{T-r})$ for $0\le r\le t\le T$. Specially, $\mcr{G}_t^0 = \sigma(\mcr{G}^{(1)}_t\cup
 \mcr{G}^{(2)}_T)$ and $\mcr{G}_T^{T-t} = \sigma(\mcr{G}^{(1)}_T\cup
 \mcr{G}^{(2)}_t)$ are two filtrations satisfying the usual hypotheses. It is
 easily seen that $\{B_t\}$ is a $(\mcr{G}_t^0)$-Brownian motion and
 $M(dt,du)$ is a $(\mcr{G}_t^0)$-Poisson random measure. Define the
 Gaussian white noise $\{W^T(dt,du): 0\le t\le T, u \in E\}$ by
 $$
 W^T([T-t, T ]\times A) = W([0, t]\times A), \qquad 0\le t\le T, A \in \mcr{B}(E).
 $$
 For $i=0,1$, define the Poisson random measures $\{N^T_i(dt,du):0\le t\le
 T, u \in U_i\}$ by:
 $$
 N^T_i([T-t, T ]\times B)= N_i ( [0, t]\times B), \qquad 0\le t\le T, B \in \mcr{B}(U_i ).
 $$
 Roughly speaking, we can consider $W^T (dt, du)$ and $N^T_i(dt,du)$ as
 the time reversal of $W(dt,du)$ and $N_i(dt,du)$, respectively.

 A real process $\{\xi_s\}_{0\le s\le T }$ is said to be
 \textit{$(\mcr{G}_t^r)$-progressive} if for any $0\le r\le t\le T$, the
 mapping $(s, \omega)\mapsto \xi_s(\omega)$ restricted to
 $[r,t]\times\Omega$ is $\mcr{B}([r,t])\times\mcr{G}_t^r$-measurable. A
 two-parameter real process $\{\zeta_s(u)\}_{0\le s\le T, u \in E}$ is said
 to be \textit{$(\mcr{G}_t^r)$-progressive} if for any $0\le r\le t\le T$,
 the restriction of $(s, u, \omega)\mapsto \zeta_s(u,\omega)$ to
 $[r,t]\times E\times \Omega$ is $\mcr{B}([r,t])\times \mcr{B}(E) \times
 \mcr{G}_t^r$-measurable.

 Let $\mcr{P}$ denote the $\sigma$-algebra on $\Omega\times [0, T]$
 generated by all real-valued left continuous processes which are
 $(\mcr{G}_t^r)$-progressive. A process $ (\xi_s)_{0\le s\le T}$ is said
 to be \textit{predictable} if the mapping $(\omega, s)\mapsto
 \xi_s(\omega)$ is $\mcr{P}$-measurable. Also a two-parameter process $\{
 \zeta_s(u)\}_{0\le s\le T, u \in E}$ is said to be \textit{predictable} if
 the mapping $(\omega,s,x)\mapsto \zeta_s(\omega,x)$ is $(\mcr{P}\times
 \mcr{B}(E))$-measurable. For the theory of time-space stochastic integrals of
 predictable two parameter processes with respect to point processes or
 random measures, readers can refer to Section II.3 in Ikeda and Watanabe (1989).
 The stochastic integrals with respect to martingale
 measures were discussed in Section 7.3 of Li (2011).
 We make the convention that the stochastic integral of a progressive
 process refers to a predictable version of the integrand.  The existence of such
 a version was briefly discussed in Section~2 of He et al. (2014).
 To simplify the following statements, we introduce several Banach spaces:
 \begin{enumerate}
  \item[(1)] $\mbb{S}_{\mcr{G},T}^{2} := \{(\xi_s)_{0\le s\le
      T}: \xi_s \mbox{ is
      }(\mcr{G}_t^r)\mbox{-progressive and }
      \|\xi\|_{\mbb{S}_T^{2}}<\infty\}$, where
      $$
      \|\xi\|_{\mbb{S}_T^{2}}^{2}=\mbf{E}\Big[\sup_{s\in[0,T]}\|\xi_s\|^{2}\Big].
      $$

 \item[(2)] $\mcr{L}_{\mcr{G},T}^{2} := \left\{(\beta_s)_{0\le
     s\le T}: \beta_s \mbox{ is
     }(\mcr{G}_t^r)\mbox{-progressive and }
     \|\beta\|_{\mcr{L}_T^{2}}<\infty \right\}$, where
      $$
      \|\beta\|_{\mcr{L}_T^{2}}^{2}=\mbf{E}\bigg\{\int_0^T\|\beta_s\|^{2}ds
      \bigg\}.
      $$

 \item[(3)] $\mcr{L}_{\mcr{G},T}^{2}(E) := \Big\{\{
     \sigma(s,u)\}_{0\le s\le T, u \in E}: \sigma(s,u) \mbox{ is } (\mcr{G}_t^r)\mbox{-progressive and }
     \|\sigma\|_{\mcr{L} _T^{2}(E)}<\infty\Big\}$, where
     $$
     \|\sigma\|_{\mcr{L}_T^{2}(E)}^{2}=\mbf{E}\left\{\int_0^T\|\sigma(s,\cdot)\|_{\mcr{L}^{2}(E)}^{2}
     ds\right\}=\mbf{E}\bigg\{\int_0^T
     ds\int_E\mathbf{T}(\sigma^{\rm T}(s,u)\sigma(s,u))\pi(du)\bigg\}.
     $$

 \item[(4)] $\mcr{L}_{\mcr{G},T}^{2}(U_0):=\Big\{\{g(s,u)\}_{0\le s\le T, u \in U_0}:g(s,u) \mbox{ is }(\mcr{G}_t^r)\mbox{-progressive and }\|g\|_{\mcr{L}_T^{2}
(U_0)}<\infty \Big\}$, where
     $$\|g\|_{\mcr{L}_T^{2}(U_0)}^{2}=\mbf{E}\left\{\int_0^T\|g(s,\cdot)\|_{\mcr{L}^{2}(U_0)}^{2}
     ds\right\}=\mbf{E}\left\{\int_0^T
     ds\int_{U_0}\|g(s,u)\|^{2}\mu_0(du)\right\}.$$

 \item[(5)] $\mcr{L}_{\mcr{G},T}^{2}(U_1):=\Big\{\{g(s,u)\}_{0\le s\le T, u \in U_1}:g(s,u)\mbox{ is }(\mcr{G}_t^r)\mbox{-progressive and } \|g\|_{\mcr{L}_T^{2}
(U_1)}<\infty \Big\}$,
     where
     $$\|g\|_{\mcr{L}_T^{2}(U_1)}^{2}=\mbf{E}\left\{\int_0^T\|g(s,\cdot)\|_{\mcr{L}^{2}(U_1)}^{2}
     ds\right\}=\mbf{E}\left\{\int_0^T
     ds\int_{U_1}\|g(s,u)\|^{2}\mu_1(du)\right\}.$$

 \item[(6)] $\mcr{L}_{\mcr{G},T}^{2}(F):=\Big\{\zeta_s(u):\ \zeta_s(u) \mbox{ is }(\mcr{G}_t^r)\mbox{-progressive and }\|\zeta\|_{\mcr{L}_T^
{2}(F)}<\infty \Big\}$,
     where
     $$\|\zeta\|_{\mcr{L}_T^{2}(F)}^{2}=\mbf{E}\left\{\int_0^T\|\zeta_s\|_{\mcr{L}^{2}(F)}^{2}
     ds\right\}=\mbf{E}\left\{\int_0^T
     ds\int_F\|\zeta_s(u)\|^{2}\nu(du)\right\}.$$

 \end{enumerate}

 Before giving the main results, we extend It\^{o} formula to the general case.
  Let $X_t$ be a $m$-dimensional stochastic process defined by:
 \beqlb\label{3.6}
 X_t\ar=\ar X_T+\int_t^T b(s)ds+\int_t^T\int_E a(s,u)W^T(\overleftarrow{ds},du) + \int_{t-}^{T-}\!\!\int_{U_0} \gamma_0(s,u)
 \tilde{N}_0^T(\overleftarrow{ds},du)\cr
 \ar\ar+ \int_{t-}^{T-}\!\!\int_{U_1}
 \gamma_1(s,u) N_1^T(\overleftarrow{ds},du) - \int_t^TZ_s dB(s)- \int_t^T\int_F\zeta_s(u)\tilde{M}(ds,du),\quad
 \eeqlb
 where $b(s)$, $a(s,u)$, $\gamma_0(s,u)$, $\gamma_1(s,u)$ and $\zeta_s(u)$ are $m$-dimensional $(\mcr{G}_t^r)$-progressive processes, $a(s,u)$ and $Z_s$
 are $(\mcr{G}_t^r)$-progressive $(m\times n)$-matrix-valued processes.

 \begin{proposition}\label{P2.1}
 For any $ f\in C^2(\mathbb{R}^m,\mathbb{R})$, we have
 \beqnn
 f(X_t) \ar=\ar f(X_T) + \int_t^T D f(X_s)b_i(s)ds+\int_t^T \int_E D f(X_s)a(s,u)W^T(\overleftarrow{ds},du) \cr
 \ar\ar
 + \int_t^T ds\int_E \mathbf{T}[a^{\rm T}(s,u) D^2 f(X_s)a(s,u)]
 \pi(du)\cr
 \ar\ar
  +\int_t^T\int_{U_0}\left[f(X_s+\gamma_0(s,u))-f(X_s)\right]\tilde{N}_0^T(\overleftarrow{ds},du)\cr
  \ar\ar+\int_t^Tds\int_{U_0}\left[f(X_s+\gamma_0(s,u))-f(X_s)-D f(X_s)\gamma_{0}(s,u)\right]\mu_0(du)\cr
 \ar\ar
  +\int_t^T\int_{U_1}\left[f(X_s+\gamma_1(s,u))-f(X_s)\right]N_1^T(\overleftarrow{ds},du)\cr
  \ar\ar - \int_t^T D f(X_s)Z_sdB(s) -\frac{1}{2}\int_t^T{\rm Tr}(Z^{\rm T}(s)D^2f(X_s)Z_s)ds\cr
  \ar\ar-\int_t^T\int_F\left[f(X_s+\zeta_s(u))-f(X_s)\right]\tilde{M}(ds,du)\cr
  \ar\ar-\int_t^Tds\int_F\left[f(X_s+\zeta_s(u))-f(X_s)-D f(X_s)\zeta_s(u)\right]\nu(du).
 \eeqnn
 \end{proposition}

  \begin{remark}\label{R2.1}
  As in He et al. (2014), we make the convention that the stochastic integral of a progressive process always refers to that of a predictable version of the integrand. Here we emphasis that the integrals in (\ref{2.1}) denote by $W^T(\overleftarrow{ds},du)$, $\tilde{N}_0^T(\overleftarrow{ds},du)$ and $N_1^T(\overleftarrow{ds},du)$ are backward ones, which can be defined as the time-reversal of the corresponding forward stochastic integrals; see He et al. (2014) for more precise explanations. Of course, the integrals with respect to $dB(s)$ and $\tilde{M}(ds,du)$ in (\ref{2.1}) are forward ones.
 \end{remark}

   Now let us introduce the backward doubly stochastic integral equation to work with. Suppose that we have the following measurable mappings:
 \beqnn
 \beta\ar :\ar [0,T]\times \mathbb{R}^m\times\mathbb{R}^{m\times n}\times\mcr{L}_{\mcr{G},T}^2(F)\mapsto \mathbb{R}^m;\cr
 \sigma\ar :\ar [0,T]\times \mathbb{R}^m\times\mathbb{R}^{m\times n}\times\mcr{L}_{\mcr{G},T}^2(F)\times E\mapsto \mathbb{R}^{m\times n};\cr
  g_0\ar :\ar [0,T]\times \mathbb{R}^m\times\mathbb{R}^{m\times n}\times\mcr{L}_{\mcr{G},T}^2(F)\times U_0\mapsto \mathbb{R}^m;\cr
  g_1\ar :\ar [0,T]\times \mathbb{R}^m\times\mathbb{R}^{m\times n}\times\mcr{L}_{\mcr{G},T}^2(F)\times U_1\mapsto \mathbb{R}^m.
 \eeqnn
  Given $Y_T\in \mcr{G}_T^0$, we consider the equation:
 \beqlb\label{2.1}
 Y_t\ar=\ar Y_T+\int_t^T\beta(s,Y_s,Z_s,\zeta_s)ds  +
 \int_t^T\int_E\sigma(s,Y_s,Z_s,\zeta_s,u)W^T(\overleftarrow{ds},du)\cr
 \ar\ar
 + \int_{t-}^{T-}\!\!\int_{U_0} g_0(s,Y_s,Z_s,\zeta_s,u)
 \tilde{N}_0^T(\overleftarrow{ds},du) + \int_{t-}^{T-}\!\!\int_{U_1}
 g_1(s,Y_s,Z_s,\zeta_s,u) N_1^T(\overleftarrow{ds},du) \cr
 \ar\ar - \int_t^TZ_s dB_s  - \int_t^T\int_F\zeta_s(u)\tilde{M}(ds,du).\quad
 \eeqlb

 \begin{definition}
 We call the process $(Y_t, Z_t,\zeta_t(u))_{0\le t\le T}$ a {\rm
 solution} to (\ref{2.1}) if it is $(\mcr{G}_t^r)$-progressive and for any
 $0\le r\le t\le T$ the equation (\ref{2.1}) is satisfied almost surely.

 \end{definition}

 \begin{condition}\label{C1}
  There exist constants $C>0$ and
 $0<\alpha<1$ such that for any $s\in [0,T]$ and $ (x_i,y_i,z_i,\zeta_i)\in \mathbb{R}^m\times\mathbb{R}^m\times\mathbb{R}^{m\times n} \times
 \mcr{L}_{\mcr{G},T}^2(F)$ with $i=1,2$,
 \beqlb\label{2.2}
 \|\beta(s,y_1,z_1,\zeta_1)-\beta(s,y_2,z_2,\zeta_2)\|^2
 \le
 C\big(\|y_1-y_2\|^2+\|z_1-z_2\|^2+\|\zeta_1-\zeta_2\|_{\mcr{L}^2(F)}^2\big)
 \eeqlb
 and
 \beqlb\label{2.3}
 \ar\ar\|\sigma(s,y_1,z_1,\zeta_1,\cdot)-\sigma(s,y_2,z_2,\zeta_2,\cdot)\|_{\mcr{L}^{2}(E)}^2
 \cr
 \ar\ar\quad+\|g_0(s,y_1,z_1,\zeta_1,\cdot)-g_0(s,y_2,z_2,\zeta_2,\cdot)\|_{\mcr{L}^{2}(U_0)}^2 \cr
 \ar\ar
 \quad+\|g_1(s,y_1,z_1,\zeta_1,\cdot)-g_1(s,y_2,z_2,\zeta_2,\cdot)\|^2_{\mcr{L}^{2}(U_1)}
 \cr
 \ar\ar \leq
 C\|y_1-y_2\|^2+\alpha\|z_1-z_2\|^2+\alpha\|\zeta_1-\zeta_2\|_{\mcr{L}^2(F)}^2.\qquad
 \eeqlb
 \end{condition}

 \begin{theorem}\label{T1}
 Suppose Condition~\ref{C1} holds. If $(Y^{(1)}_t,Z^{(1)}_t,\zeta^{(1)}_t(u))$
 and $(Y^{(2)}_t,Z^{(2)}_t,\zeta^{(2)}_t(u))$ are solutions to (\ref{2.1}) with $Y^{(1)}_T
 = Y^{(2)}_T$ a.s., then
 \beqlb\label{2.4a}
 \mbf{P}\left(Y_t^{(1)}=Y_t^{(2)} \mbox{ for all } t\in [0,T]\right) = 1
 \eeqlb
 and
 \beqlb\label{2.4}
 \|Z^{(1)}-Z^{(2)}\|_{\mcr{L}_T^{2}} +
 \|\zeta^{(1)}-\zeta^{(2)}\|_{\mcr{L}_T^{2}(F)} = 0.
 \eeqlb
 \end{theorem}
 \proof  Let $(\bar{Y}_t, \bar{Z}_t,
 \bar{\zeta}_t(u)) = (Y^{(1)}_t-Y^{(2)}_t, Z^{(1)}_t-Z^{(2)}_t, \zeta^{(1)}_t(u) -
 \zeta^{(2)}_t(u))$. From (\ref{2.1}) we get
 \beqlb\label{2.6a}
 \bar{Y}_t
 \ar=\ar
 \int_t^T\bar{\beta}(s)ds + \int_t^T\int_E\bar{\sigma}(s,u)
 W^T(\overleftarrow{ds},du)
 +\int_{t-}^{T-}\int_{U_0} \bar{g_0}(s,u) \tilde{N}_0^T(\overleftarrow{ds},du)\cr
 \ar\ar +
 \int_{t-}^{T-}\int_{U_1} \bar{g_1}(s,u) N_1^T(\overleftarrow{ds},du)
 - \int_t^T\bar{Z}_s dB_s -
 \int_t^T\int_F\bar{\zeta}_{s}(u)\tilde{M}(ds,du),
 \eeqlb
 where
 \beqnn
 \bar{\beta}(s) \ar=\ar
 \beta(s,Y^{(1)}_t,Z^{(1)}_t,\zeta^{(1)}_t)-\beta(s,Y^{(2)}_t,Z^{(2)}_t,\zeta^{(2)}_t),\cr
 \bar{\sigma}(s,u) \ar=\ar
 \sigma(s,Y^{(1)}_t,Z^{(1)}_t,\zeta^{(1)}_t,u)-\sigma(s,Y^{(2)}_t,Z^{(2)}_t,\zeta^{(2)}_t,u),\cr
 \bar{g_0}(s,u) \ar=\ar
 g_0(s,Y^{(1)}_t,Z^{(1)}_t,\zeta^{(1)}_t,u)-g_0(s,Y^{(2)}_t,Z^{(2)}_t,\zeta^{(2)}_t,u),\cr
 \bar{g_1}(s,u) \ar=\ar
 g_1(s,Y^{(1)}_t,Z^{(1)}_t,\zeta^{(1)}_t,u)-g_1(s,Y^{(2)}_t,Z^{(2)}_t,\zeta^{(2)}_t,u).
 \eeqnn
 By Proposition~\ref{P2.1}, we
 have
 \beqnn
 \|\bar{Y}_t\|^2
 \ar=\ar
 2\int_t^T\langle\bar{Y}_s,\bar{\beta}(s)\rangle ds +
 2\int_t^T\int_E\langle\bar{Y}_s,\bar{\sigma}(s,u)\rangle W^T(\overleftarrow{ds},du)+
 \int_t^T \|\bar{\sigma}(s,\cdot)\|_{\mcr{L}^{2}(E)}^2ds\cr
 \ar\ar
 +\int_t^T\int_{U_0}[2\langle\bar{Y}_{s}, \bar{g_0}(s,u)\rangle + \|
 \bar{g_0}(s,u)\|^2]\tilde{N}_0^T(\overleftarrow{ds},du) +\int_t^T\|\bar{g_0}(s,\cdot)\|_{\mcr{L}^{2}(U_0)}^2 ds\cr
 \ar\ar
 +2\int_t^T\int_{U_1}\langle\bar{Y}_{s},
 \bar{g_1}(s,u)\rangle N_1^T(\overleftarrow{ds},du)+\int_t^T\int_{U_1}
 \|\bar{g_1}(s,u)\|^2N_1^T(\overleftarrow{ds},du)\cr
 \ar\ar
 -\, 2\int_t^T\langle\bar{Y}_s,\bar{Z}_s\rangle dB_s - \int_t^T\|\bar{Z}_s\|^2 ds - \int_t^T\|\bar{\zeta}_s\|_{\mcr{L}^{2}(F)}^2
 ds\cr
 \ar\ar
 -\int_t^T\int_F[2\langle\bar{Y}_{s},\bar{\zeta}_{s}(u)\rangle+\|\bar{\zeta}_{s}(u)\|^2]\tilde{M}(ds,du).
 \eeqnn
 From Cauchy's inequality, for any $a,b >0$ we have
 \beqnn
 \lefteqn{\mbf{E}\big[\|\bar{Y}_t\|^2\big] + \mbf{E}\Big[\int_t^T \|\bar{Z}_s\|^2
 ds\Big] + \mbf{E}\Big[\int_t^T \|\bar{\zeta}_s\|_{\mcr{L}^2(F)}^2 ds\Big]}\qquad\ar\ar\cr
 \ar=\ar
 \mbf{E}\left\{2\int_t^T\langle \bar{Y}_s,\bar{\beta}(s)\rangle ds\right\}+
 \mbf{E}\left\{\int_t^T \|\bar{\sigma}(s,\cdot)\|_{\mcr{L}^{2}(E)}^2 ds\right\}+ \mbf{E}\left\{\int_t^T
 \|\bar{g_0}(s,\cdot)\|_{\mcr{L}^{2}(U_0)}^2ds\right\}\cr
 \ar\ar
 +
 \mbf{E}\left\{2\int_t^Tds\int_{U_1}\langle \bar{Y}_s,
 \bar{g_1}(s,u)\rangle \mu_1(du)\right\}
 + \mbf{E}\left\{\int_t^T\|\bar{g_1}(s,\cdot)\|_{\mcr{L}^{2}(U_1)}^2 ds
 \right\}\cr
 \ar\le\ar
 (\frac{1}{a}+\frac{\mu_1(U_1)}{b})\mbf{E}\left\{\int_t^T\|\bar{Y}_s\|^2ds\right\}
 + a\mbf{E}\left\{\int_t^T\|\bar{\beta}(s)\|^2ds\right\}
 + \mbf{E}\left\{\int_t^T\|\bar{\sigma}(s,\cdot)\|_{\mcr{L}^{2}(E)}^2  ds\right\}\cr
 \ar\ar
 +\, \mbf{E}\left\{\int_t^T \|\bar{g_0}(s,\cdot)\|_{\mcr{L}^{2}(U_0)}^2ds\right\}+ (1+b)\mbf{E}\left\{\int_t^T \|\bar{g_1}(s,\cdot)\|_{\mcr{L}^{2}(U_1)}^2 ds\right\}.
 \eeqnn
 Since $\mu_1$ is a finite measure, by H\"{o}lder's inequality and
 Condition~\ref{C1},
  \beqnn
 \lefteqn{\mbf{E}\big[\|\bar{Y}_t\|^2\big] + \mbf{E}\Big[\int_t^T \|\bar{Z}_s\|^2
 ds\Big] + \mbf{E}\Big[\int_t^T \|\bar{\zeta}_s\|_{\mcr{L}^2(F)}^2 ds
 \Big]}\ar\ar\cr
 \ar\ar
 \le
 \Big(\frac{1}{a}+\frac{\mu_1(U_1)}{b}\Big)\mbf{E}\left\{\int_t^T\|\bar{Y}_s\|^2ds\right\}
 + Ca\mbf{E}\left\{\int_t^T[\|\bar{Y}_s\|^2 +
 \|\bar{Z}_s\|^2+\|\bar{\zeta}_s\|_{\mcr{L}^2(F)}^2] ds\right\}\cr
 \ar\ar\quad
  +\, (1+b)\mbf{E}\left\{\int_t^T [C\|\bar{Y}_s\|^2+\alpha\|\bar{Z}_s\|^2 +
 \alpha\|\bar{\zeta}_s\|_{\mcr{L}^2(F)}^2]ds\right\}.
 \eeqnn
 Here we can choose $a,b$ small enough such that $\hat{\alpha} := Ca +
 \alpha + b\alpha< 1$. Then
 \beqnn
 \mathbf{E}\big[\|\bar{Y}_t\|^2\big]+(1-\hat{\alpha})\mbf{E}\left\{\int_t^T\|\bar{Z}_s\|^2
 ds\right\}+(1-\hat{\alpha})\mbf{E}\left\{\int_t^T
 \|\bar{\zeta}_s\|_{\mcr{L}^2(F)}^2 ds\right\}\cr
  \leq \big[1/a+1/b+C(1+a+b)\big]\mbf{E}\left\{\int_t^T\|\bar{Y}_s\|^2ds\right\}.
 \eeqnn
 By Gronwall's lemma, we have
 \beqnn
 \mbf{E}\big[\|\bar{Y}_t\|^2\big] + \mbf{E}\left\{\int_t^T\|\bar{Z}_s\|^2 ds\right\} +
 \mbf{E}\left\{\int_t^T \|\bar{\zeta}_s\|_{\mcr{L}^2(F)}^2ds\right\} =
 0.
 \eeqnn
 This implies (\ref{2.4}). Then for any fixed $t\in [0,T]$, the six terms
 on the right-hand side of (\ref{2.6a}) vanish almost surely. Since each of
 the six terms is right-continuous or left-continuous, they almost surely
 vanish for all $t\in [0,T]$. We have finished the proof.
 \qed

 \section{Existence}
 \setcounter{equation}{0}
 \medskip

 In this section, we study the existence of solutions to (\ref{2.1}). For any $0\le r\le t\le T$ we define the natural $\sigma$-algebras:
 \beqnn
 \mcr{F}^{BM}_{r,t} \ar=\ar
 \sigma(\{B(s)-B(r), M((r,s]\times A): r\le s\le t, A\in \mbf{B}(F)\})\vee
 \mcr{N}, \cr
 \mcr{F}^{WN}_{r,t} \ar=\ar
 \sigma(\{W((r,s]\times A), N_0((r,s]\times B), N_1((r,s]\times C): \cr
 \ar\ar\qquad
 r\le s\le t, A\in\mcr{B}(E), B\in \mcr{B}(U_0), C\in \mcr{B}(U_1)\})\vee
 \mcr{N},
 \eeqnn
 where $\mcr{N}$ denotes the totality of $\mbf{P}$-null sets. For
 simplicity, we write $\mcr{F}^{BM}_t = \mcr{F}^{BM}_{0,t}$ and
 $\mcr{F}^{WN}_t = \mcr{F}^{WN}_{0,t}$. Let
 $\mcr{F}_t^r = \mcr{F}^{BM}_t \vee \mcr{F}^{WN}_{T-r}$ for $0\le r\le t\le T$.
 Similarly, we can define $\mathbb{S}_{\mcr{F},T}^2$, $\mcr{L}_{\mcr{F},T}^2$,
 $\mcr{L}_{\mcr{F},T}^2(E)$, $\mcr{L}_{\mcr{F},T}^2(U_0)$,
 $\mcr{L}_{\mcr{F},T}^2(U_1)$, $\mcr{L}_{\mcr{F},T}^2(F)$ like those in the last section but with $\{\mcr{G}_t^r: 0\leq r\leq t\leq T\}$ replaced by
 $\{\mcr{F}_t^r: 0\leq r\leq t\leq T\}$.

 \begin{theorem}\label{T2}
 Suppose Condition~\ref{C1} holds. Then there exists a solution $(Y_t,Z_t,\zeta_t(u))$ to
 (\ref{2.1}) in $\mathbb{S}_{\mcr{F},T}^2\times \mcr{L}_{\mcr{F},T}^2\times
 \mcr{L}_{\mcr{F},T}^2(F)$.
 \end{theorem}

 Obviously, combining this theorem with Theorem~\ref{T1}, we have solution to
 (\ref{2.1}) exists uniquely in $\mathbb{S}_{\mcr{G},T}^2\times \mcr{L}_{\mcr{G},T}^2\times
 \mcr{L}_{\mcr{G},T}^2(F)$.
 Before giving the proof of Theorem~\ref{T2}, we introduce a lemma about
 the solution to some simple backward doubly stochastic equation, which is
 very important in the proof of this theorem.

 \begin{lemma}\label{L1}
 Let $\beta\in \mcr{L}_{\mcr{F},T}^1$, $\sigma\in
 \mcr{L}_{\mcr{F},T}^2(E)$, $g_0\in \mcr{L}_{\mcr{F},T}^2(U_0)$ and $g_1\in
 \mcr{L}_{\mcr{F},T}^2(U_1)$.
 Then for any $Y_T\in \mcr{F}_{T}^0$ with finite second moment, there
 exists a unique solution $(Y_t,Z_t,\zeta_t(u))\in \mathbb{S}_{\mcr{F},T}^2
 \times \mcr{L}_{\mcr{F},T}^2\times \mcr{L}_{\mcr{F},T}^2(F)$ to the
 following equation:
 \beqlb\label{3.4}
 Y_t
 \ar=\ar
 Y_T + \int_t^T\beta(s)ds + \int_t^T\int_E \sigma(s,u)
 W^T(\overleftarrow{ds},du) + \int_{t-}^{T-}\int_{U_0}
 g_0(s,u)\tilde{N}_0^T(\overleftarrow{ds},du)\cr
 \ar\ar
 +\int_{t-}^{T-}\int_{U_1} g_1(s,u)N_1^T(\overleftarrow{ds},du) -
 \int_t^TZ_s dB_s-\int_t^T\int_F\zeta_s(u)\tilde{M}(ds,du).
 \eeqlb
 \end{lemma}

 \proof The uniqueness of the solution follows from Theorem~\ref{T1}. Recall
  $\mcr{F}^0_t = \mcr{F}^{BM}_t\vee \mcr{F}^{WN}_{T}$ for $0\le t\le
 T$. Observe that
 \beqlb\label{3.9}
 \Psi_T\ar:=\ar Y_T + \int_0^T\int_E\sigma(s,u)W^T(\overleftarrow{ds},du) +
 \int_{0-}^{T-}\int_{U_0} g_0(s,u)\tilde{N}_0^T(\overleftarrow{ds},du)\cr
 \ar\ar\qquad
 + \int_0^T\beta(s)ds + \int_{0-}^{T-}\int_{U_1}
 g_1(s,u)N_1^T(\overleftarrow{ds},du)
 \eeqlb
 is $\mcr{F}^0_T$-measurable. Then we can define a Doob's martingale:
 \beqnn
 M_t=\mbf{E}[\Psi_T|\mcr{F}^0_t], \qquad 0\le t\le T.
 \eeqnn
 Since $\mcr{F}_t^t\subset\mcr{F}^0_t$, from (\ref{3.9}) we have
 \beqlb\label{3.10}
 M_t\ar=\ar Y_t + \int_0^t\beta(s)ds + \int_0^t\int_E \sigma(s,u) W^T(\overleftarrow{ds},du)\cr
 \ar\ar +
 \int_{0-}^{t-}\int_{U_0} g_0(s,u)
 \tilde{N}_0^{T}(\overleftarrow{ds},du)
 + \int_{0-}^{t-}\int_{U_1} g_1(s,u)
 N_1^T(\overleftarrow{ds},du),
 \eeqlb
 where $Y_t=\mbf{E}[\Xi(t)|\mcr{F}^0_t]$ and
 \beqlb\label{3.10a}
 \Xi(t)\ar=\ar Y_T + \int_t^T\beta(s)ds + \int_t^{T}\int_E\sigma(s,u) W^T(\overleftarrow{ds},du)
 \cr
 \ar\ar + \int_{t-}^{T-}\int_{U_0} g_0(s,u)
 \tilde{N}_0^{T}(\overleftarrow{ds},du)
 + \int_{t-}^{T-}\int_{U_1}
 g_1(s,u)N_1^T(\overleftarrow{ds},du).
 \eeqlb
 By the martingale representation theorem, see Lemma 2.3 in Tang and Li (1994), there exist
 $(\mcr{F}^0_t)$-progressive processes $\{Z_s\}$ and $\{\zeta_s(u)\}$ such
 that
 \beqnn
 M_t= M_0+\int_0^t Z_s dB_s+ \int_0^t \int_F \zeta_s(du)\tilde{M}(ds,du)
 \eeqnn
 and hence
 \beqlb\label{3.11}
 M_T = M_t + \int_t^T Z_s dB_s + \int_t^T \int_F
 \zeta_s(u)\tilde{M}(ds,du).
 \eeqlb
Since $M_T=\Psi_T$, we can substitute (\ref{3.9}) and (\ref{3.10}) into
 (\ref{3.11}) to obtain (\ref{3.4}). Finally, we need to prove for any
 $0\le r\le T$ the process $(Y_t,Z_t,\zeta_t(u))_{r\le t\le T, u\in F}$ is
 $(\mcr{F}_t^r)$-progressive. Observe that
  $$
 Y_r = \mbf{E}[\Xi(r)|\mcr{F}^0_r]
 =
 \mbf{E}[\Xi(r)|\mcr{F}^{BM}_r\vee \mcr{F}^{WN}_{T}]
 =
\mbf{E}[\Xi(r)|\mcr{F}_r^r\vee \mcr{F}^{WN}_{T-r,T}],
 $$
 where $\mcr{F}^{WN}_{T-r}$ and $ \mcr{F}^{WN}_{T-r,T}$ are independent.
 By (\ref{3.10a}) it is easy to see that $\Xi(r)$ is independent of
 $\mcr{F}^{WN}_{T-r,T}$. Then we have $Y_r = \mbf{E}[\Xi(r)|\mcr{F}_r^r]$,
 which is $\mcr{F}_r^r$-measurable. By (\ref{3.4}) we have
  \beqnn
 \int_r^T Z_s dB_s+\int_r^T \int_F \zeta_s(u)
 \tilde{M}(ds,du)
 \ar=\ar \int_r^T\int_E \sigma(s,u) W^T(\overleftarrow{ds},du) +
 \int_{r-}^{T-}\int_{U_0} g_0(s,u) \tilde{N}_0^{T}(\overleftarrow{ds},du)
 \cr
 \ar\ar
 +Y_T - Y_r + \int_r^T\beta(s)ds + \int_{r-}^{T-}\int_{U_1} g_1(s,u)
 N_1^T(\overleftarrow{ds},du).
 \eeqnn
 Then by the uniqueness of the martingale representation, the process
 $(Z_t,\zeta_t(u))$ has an $(\mcr{F}^r_t)$-progressive version. Since each
 term in (\ref{3.4}) is right or left continuous, the process $(Y_t)$ is
 $(\mcr{F}^r_t)$-progressive.
 \qed

 \noindent\textit{Proof of Theorem~\ref{T2}.~} \ We shall use a Picard
 iteration argument to construct a solution to (\ref{2.1}). Let $Y^{(0)}_t =
 Z^{(0)}_t = \zeta^{(0)}_t(u)\equiv 0$. By Lemma~\ref{L1}, for any $n\geq
 0$ there exists a unique solution $(Y^{(n+1)}_t, Z^{(n+1)}_t,
 \zeta^{(n+1)}_t(u))$ to the following equation:
 \beqnn
 Y^{(n+1)}_t\ar=\ar Y_T + \int_t^T
 \beta(s,Y^{(n)}_s,Z^{(n)}_s,\zeta^{(n)}_s) ds + \int_t^T\int_E
 \sigma(s,Y^{(n)}_s,Z^{(n)}_s,\zeta^{(n)}_s,u)
 W^T(\overleftarrow{ds},du)\cr
 \ar\ar
 + \int_{t-}^{T-}\int_{U_0} g_0(s,Y^{(n)}_{s},Z^{(n)}_{s},\zeta^{(n)}_{s},u) \tilde{N}_0^T(\overleftarrow{ds},du) - \int_t^T Z^{(n+1)}_s dB_s \cr
 \ar\ar
 + \int_{t-}^{T-}\int_{U_1} g_1(s,Y^{(n)}_{s},Z^{(n)}_{s},\zeta^{(n)}_{s},u) N_1^T(\overleftarrow{ds},du)
 - \int_t^T\int_F\zeta^{(n+1)}_{s}(u)\tilde{M}(ds,du).
 \eeqnn
 Let $\bar{Y}^{(n+1)}_t=Y^{(n+1)}_t-Y^{(n)}_t$, $\bar{Z}^{(n+1)}_t =
 Z^{(n+1)}_t-Z^{(n)}_t$ and $\bar{\zeta}^{(n+1)}_t(u) = \zeta^{(n+1)}_t(u)
 - \zeta^{(n)}_t(u)$. From (\ref{2.1}) we have
 \beqnn
 \bar{Y}^{(n+1)}_t\ar=\ar\int_t^T\bar{\beta}^{(n)}(s)ds
 + \int_t^T\int_E \bar{\sigma}^{(n)}(s,u) W^T(\overleftarrow{ds},du) + \int_{t-}^{T-}\int_{U_0} \bar{g_0}^{(n)}(s,u) \tilde{N}_0^T(\overleftarrow{ds},du)\cr
 \ar\ar
+\int_{t-}^{T-}\int_{U_1} \bar{g_1}^{(n)}(s,u) N_1^T(\overleftarrow
{ds},du)
 - \int_t^T \bar{Z}^{(n+1)}_s dB_s - \int_t^T\int_F \bar{\zeta}^{(n+1)}_{s}(u) \tilde{M}(ds,du),
 \eeqnn
 where
 \beqnn
 \bar{\beta}^{(n)}(s)\ar=\ar
 \beta(s,Y^{(n)}_s,Z^{(n)}_s,\zeta^{(n)}_s) - \beta(s,Y^{(n-1)}_s,Z^{(n-1)}_s,\zeta^{(n-1)}_s),\cr
 \bar{\sigma}^{(n)}(s,u)\ar=\ar
 \sigma(s,Y^{(n)}_s,Z^{(n)}_s,\zeta^{(n)}_s,u) - \sigma(s,Y^{(n-1)}_s,Z^{(n-1)}_s,\zeta^{(n-1)}_s,u),\cr
 \bar{g_0}^{(n)}(s,u)\ar=\ar
 g_0(s,Y^{(n)}_s,Z^{(n)}_s,\zeta^{(n)}_s,u) - g_0(s,Y^{(n-1)}_s,Z^{(n-1)}_s,\zeta^{(n-1)}_s,u),\cr
 \bar{g_1}^{(n)}(s,u)\ar=\ar
 g_1(s,Y^{(n)}_s,Z^{(n)}_s,\zeta^{(n)}_s,u) - g_1(s,Y^{(n-1)}_s,Z^{(n-1)}_s,\zeta^{(n-1)}_s,u).
 \eeqnn
 According to Proposition~\ref{P2.1},, we have
 \beqnn
 \|\bar{Y}^{(n+1)}_t\|^2 \ar=\ar
 2\int_t^T\langle\bar{Y}^{(n+1)}_s,\bar{\beta}^{(n)}(s)\rangle ds + 2\int_t^T\int_E \langle\bar{Y}^{(n+1)}_s, \bar{\sigma}^{(n)}(s,u)\rangle W^T(\overleftarrow{ds},du)\cr
 \ar\ar
 + \int_t^T\|\bar{\sigma}^{(n)}(s,\cdot)\|_{\mcr{L}^2(E)}^2ds+\int_t^T\|\bar{g_0}^{(n)}(s,\cdot)\|_{\mcr{L}^2(U_0)}^2ds\cr
 \ar\ar
 +\, \int_{t-}^{T-}\int_{U_0}[2\langle\bar{Y}^{(n+1)}_s, \bar{g_0}^{(n)}(s,u)\rangle+\|\bar{g_0}^
{(n)}(s,u)\|^2] \tilde{N}_0^T(\overleftarrow{ds},du)\cr
 \ar\ar
 +\, \int_{t-}^{T-}\int_{U_1}[2\langle\bar{Y}^{(n+1)}_s, \bar{g_1}^{(n)}(s,u)\rangle + \|\bar{g_1}^{(n)}
(s,u)\|^2]N_1^T(\overleftarrow{ds},du)\cr
 \ar\ar
 -\,2\int_t^T\langle\bar{Y}^{(n+1)}_s,\bar{Z}^{(n+1)}_s\rangle dB_s - \int_t^T \|\bar{Z}^{(n+1)}_s\|^2 ds- \int_t^T\|\bar{\zeta}^{(n+1)}_s\|_{\mcr{L}^2(F)}^2 ds \cr
 \ar\ar - \int_t^T\int_F [2\langle\bar{Y}^{(n+1)}_s,\bar{\zeta}^{(n+1)}_{s}(u)\rangle + \|\bar{\zeta}^{(n+1)}_{s}(u)\|^2]
\tilde{M}(ds,du).
 \eeqnn
 It follows that
 \beqnn
 \lefteqn{\mbf{E}\big[|\bar{Y}^{(n+1)}_t|^2\big]+\mbf{E}\left\{\int_t^T\|\bar{Z}^{(n+1)}_s\|^2ds\right\} + \mbf{E}\left\{\int_t^T \|\bar{\zeta}^{(n+1)}_s
\|_{\mcr{L}^2(F)}^2 ds\right\}}\qquad\ar\ar\cr
 \ar=\ar
 \mbf{E}\left\{2\int_t^T\langle\bar{Y}^{(n+1)}_s,\bar{\beta}^{(n)}(s)\rangle ds\right\} +
 \mbf{E}\left\{\int_t^T\|\bar{\sigma}^{(n)}(s,\cdot)\|_{\mcr{L}^2(E)}^2ds\right\}\cr
 \ar\ar
 +\, \mbf{E}\left\{\int_t^T\|\bar{g_0}^{(n)}(s,\cdot)\|_{\mcr{L}^2(U_0)}^2 ds\right\}+ \mbf{E}\left\{\int_{t}^{T}\|\bar{g_1}^{(n)}(s,\cdot)\|_{\mcr{L}^2(U_1)}^2 ds\right\}\cr
 \ar\ar
 +\, \mbf{E}\left\{2\int_{t}^{T}ds\int_{U_1} \langle\bar{Y}^{(n+1)}_s,
\bar{g_1}^{(n)}(s,u)\rangle \mu_1(du)\right\}.
 \eeqnn
 By integration by parts, one can see, for any $\lambda> 0$,
 \beqnn
 \lambda\int_t^T e^{\lambda s}\mbf{E}\big[\|\bar{Y}^{(n+1)}_s\|^2\big]ds\ar=\ar e^{\lambda s}\mbf{E}\big[\|\bar{Y}^{(n+1)}_s\|^2\big]\big|_t^T-\int_t^Te^{\lambda
s}d\mbf{E}[\|\bar{Y}^{(n+1)}_s\|^2]\cr
 \ar=\ar
 -e^{\lambda t}\mathbf{E}\big[\|\bar{Y}^{(n+1)}_t\|^2\big]+\mbf{E}\left\{2\int_t^Te^{\lambda s}\langle\bar{Y}^{(n+1)}_s,\bar{\beta}^{(n)}(s)\rangle ds\right\}\cr
 \ar\ar +\mbf{E}\left\{\int_t^Te^{\lambda
 s}\|\bar{\sigma}^{(n)}(s,\cdot)\|_{\mcr{L}^2(E)}^2ds\right\}\cr
 \ar\ar +
 \mbf{E}\left\{\int_t^T e^{\lambda s}\|\bar{g_0}^{(n)}(s,\cdot)\|_{\mcr{L}^2(U_0)}^2ds\right\}\cr
 \ar\ar
 +\, \mbf{E}\left\{\int_{t}^{T}e^{\lambda s}\|\bar{g_1}^{(n)}(s,\cdot)\|_{\mcr{L}^2(U_1)}^2ds\right\}\cr
 \ar\ar
 + \mbf{E}\left\{2\int_{t}^{T}e^{\lambda s}ds\int_{U_1}\langle\bar{Y}^{(n+1)}_s, \bar{g_1}^{(n)}(s,u) \rangle \mu_1(du)\right\} \cr
 \ar\ar - \mbf{E}\left\{\int_t^T\|\bar{Z}^{(n+1)}_s\|^2
e^{\lambda s}ds\right\}
 -\mbf{E}\left\{\int_t^Te^{\lambda s}\|\bar{\zeta}^{(n+1)}_s\|_{\mcr{L}^2(F)}^2ds\right\}.
  \eeqnn
 By H\"{o}lder's inequality, for any $a,b>0$ we have
 \beqnn
 \lefteqn{\lambda\int_t^T e^{\lambda s}\mbf{E}\big[\|\bar{Y}^{(n+1)}_s\|^2\big]ds + \mbf{E}\left\{\int_t^T e^{\lambda s}\|\bar{Z}^{(n+1)}_s\|^2ds\right\} +
\mbf{E}\left\{\int_t^Te^{\lambda s}\|\bar{\zeta}^{(n+1)}_s\|_{\mcr{L}^2(F)}^2 ds\right\}}\qquad\ar\ar\cr
 \ar\leq\ar
 (1/a+1/b)\int_t^Te^{\lambda
 s}\mbf{E}\big[\|\bar{Y}^{(n+1)}_s\|^2\big]ds+ \mbf{E}\left\{a\int_t^Te^{\lambda s}
 \|\bar{\beta}^{(n)}(s)\|^2ds\right\}\cr
 \ar\ar +
 \mbf{E}\bigg\{\int_t^Te^{\lambda
 s}\big[\|\bar{\sigma}^{(n)}(s,u)\|_{\mcr{L}^2(E)}^2
 +\|\bar{g_0}^{(n)}(s,u)\|_{\mcr{L}^2(U_0)}^2\big]ds\bigg\}\cr
  \ar\ar
 +(1+b)\mbf{E}\left\{\int_{t}^{T}e^{\lambda
 s}\|\bar{g_1}^{(n)}(s,u)\|_{\mcr{L}^2(U_1)}^2ds\right\}.
  \eeqnn
 Using Condition~\ref{C1}, we have
 \beqnn
 \lefteqn{\lambda\int_t^T e^{\lambda s}\mbf{E}\big[\|\bar{Y}^{(n+1)}_s\|^2\big]ds + \mbf{E}\left\{\int_t^Te^{\lambda s}\| \bar{Z}^{(n+1)}_s\|^2ds\right\} +
\mbf{E}\left\{\int_t^Te^{\lambda s}\|\bar{\zeta}^{(n+1)}_s\|_{\mcr{L}^2(F)}^2 ds\right\}}\ar\ar\cr
 \ar\leq\ar
 (1/a+1/b)\int_t^Te^{\lambda
 s}\mbf{E}\big[\|\bar{Y}^{(n+1)}_s\|^2\big]ds
 +\, aC\mbf{E}\left\{\int_t^T e^{\lambda s}\big[ \|\bar{Y}^{(n)}_s\|^2 +
 \|\bar{Z}^{(n)}_s\|^2 +
 \|\bar{\zeta}^{(n)}_s\|_{\mcr{L}^2(F)}^2\big]ds\right\}\cr
 \ar\ar
 +\, (1+b)\mbf{E}\left\{\int_t^Te^{\lambda s}\big[C\|\bar{Y}^{(n)}_s\|^2 +
 \alpha\|\bar{Z}^{(n)}_s\|^2 +
 \alpha\|\bar{\zeta}^{(n)}_s\|_{\mcr{L}^2(F)}^2\big] ds\right\}.
 \eeqnn
 Then
 \beqnn
 \lefteqn{\mbf{E}\left\{\int_t^T e^{\lambda s}
 \big[(\lambda-1/a-1/b)\|\bar{Y}^{(n+1)}_s\|^2 +\|\bar{Z}^{(n+1)}_s\|^2+ \|\bar{\zeta}^{(n+1)}_s\|_{\mcr{L}^2(F)}^2 \big]ds\right\}}\ar\ar\cr
 \ar\leq\ar
 \mbf{E}\left\{\int_t^T e^{\lambda s}\big[(a+b+1)C \|\bar{Y}^{(n)}_s\|^2 +
 (aC+b\alpha+\alpha)(\|\bar{Z}^{(n)}_s\|^2 +
 \|\bar{\zeta}^{(n)}_s\|_{\mcr{L}^2(F)}^2)\big]ds\right\}.
 \eeqnn
 Let $a,b$ be small enough such that $aC+b\alpha+\alpha<1$. Then choose $\lambda>0$ large enough such that $$
 \lambda-\frac{1}{a}-\frac{1}{b}>\frac{(a+b+1)C}{aC+b\alpha+\alpha}>0.
 $$
 It follows that
 \beqnn
 \lefteqn{\mbf{E}\left\{\int_t^T e^{\lambda s}\Big[\Big(\lambda-\frac{1}{a} - \frac{1}{b}\Big)\|\bar{Y}^{(n+1)}_s\|^2 +\|\bar{Z}^{(n+1)}_s\|^2 +
 \|\bar{\zeta}^{(n+1)}_s\|_{\mcr{L}^2(F)}^2\Big]ds\right\}}\qquad\ar\ar\cr
 \ar\leq\ar
 (aC+b\alpha+\alpha)\mbf{E}\left\{\int_t^T e^{\lambda s}\Big[\Big(\lambda - \frac{1}{a} - \frac{1}{b}\Big)\|\bar{Y}^{(n)}_s\|^2 + \|\bar{Z}^{(n)}
_s\|^2 + \|\bar{\zeta}^{(n)}_s\|_{\mcr{L}^2(F)}^2\Big]ds\right\}\cr
 \ar\leq\ar \cdots\cr
 \ar\leq\ar
 (aC+b\alpha+\alpha)^n\mbf{E}\left\{\int_t^T e^{\lambda s} \Big[\Big(\lambda-\frac{1}{a}-\frac{1}{b}\Big)\|\bar{Y}^{(1)}_s\|^2 + \|\bar{Z}^{(1)}_s\|^2 + \|\bar{\zeta}^{(1)}_s\|_{\mcr{L}^2(F)}^2\Big]ds\right\}.
 \eeqnn
 Since the right-hand side of the inequality is summable, we see that $\{(Y^{(n)}_s, Z^{(n)}_s, \zeta^{(n)}_s(u))\}$ is a Cauchy sequence. By Burkholder-Davis-Gundy Inequality, it is easy to see $Y^{(n)}_s$ is also a Cauchy sequence in $\mathbb{S}_{\mathscr{F},T}^2$.  Then
it converges in $\mathbb{S}_{\mcr{F},T}^2\times \mcr{L}_{\mcr{F},T}^2\times \mcr{L}_{\mcr{F},T}^2(F)$ to some process $(Y_s, Z_s, \zeta_s(u))$,
which is clearly a solution to (\ref{2.1}). Then we have finished the proof. \qed

\section{Comparison theorems}
\setcounter{equation}{0}
\medskip

 Comparison theorems are very important in both theory and applications.
 For instance, if you want to earn more money from a complete capital market in the future time $T$, you should either invest more money in the market at time $0$ or improve your investment policy.
 This section will mainly introduce several comparison theorems under Condition~\ref{C1}.
 There are two classical ways to prove comparison theorems in the theory of BSDEs; see Situ (2005, p.243-250).
 One is transforming the BSDE into a summation of a non-negative processes and a martingale under a new probability measure.
 Then the desired results can be gotten by taking conditional expectation under the new probability measure.
 Another one is called "a duality method" which mainly by constructing a relative forward SDE (FSDE).
 Applying It\^{o} formula to the multiplication of the solutions of these two stochastic equations (FBSDE), we will get a new process which is a summation of a non-negative processes and a martingale.
 Similarly, we get the comparison theorem by taking conditional expectation.
 Actually, both of these two methods come from the same ideas.

 Unfortunately, effected by backward integral parts in  (\ref{2.1}), neither of these two methods works.
 Here we use another method to get comparison theorems under some conditions which are not really stronger than those in BSDEs.
 The main difficulty is to deal with the influence of $\zeta_s$ to the drift coefficient and backward integrals.
 We divide the influence into several parts and deal with them one by one.
 Here we only consider the one-dimensional case, comparison theorem for multi-dimensional case is still an open problem; see Peng (1999).
 Firstly, we give a simple comparison theorem about the non-positivity of solution to the following one-dimensional BDSDEs, which can be used to derive other results.
 \beqlb\label{4.1}
 Y_t\ar=\ar Y_T+\int_t^T\beta(s,Y_s,Z_s,\zeta_s)ds +
 \int_t^T\int_E\sigma(s,Y_s,Z_s,u)W^T(\overleftarrow{ds},du)\cr
 \ar\ar
 + \int_{t-}^{T-}\!\!\int_{U_0} g_0(s,Y_{s},Z_{s},u)
 \tilde{N}_0^T(\overleftarrow{ds},du)+ \int_{t-}^{T-}\!\!\int_{U_1}
 g_1(s,Y_{s},Z_{s},u) N_1^T(\overleftarrow{ds},du)\cr
 \ar\ar
 - \int_t^TZ_s dB_s- \int_t^T\!\!\int_F\zeta_s(u)\tilde{M}(ds,du).
 \eeqlb
 \begin{lemma}\label{Comparison.L1}
 Suppose Condition~\ref{C1} holds, $(Y_t,Z_t, \zeta_t)$ is a solution to BDSDE (\ref{4.1}) and
 \begin{enumerate}

 \item[(1)] both $ y+g_0(s,y,z,u)$ and $ y+g_1(s,y,z,u)$ are non-positive for any $y\in(-\infty,0]$;

 \item[(2)] there exist some constants $C>0$ and $0<\alpha<1$ such that for any $s\in[0,T]$,
 \beqnn
 \|\sigma(s,y,z,\cdot) \|_{\mcr{L}^{2}(E)}^2+\|g_0(s,y,z,\cdot)\|_{\mcr{L}^{2}(U_0)}^2 +\|g_1(s,y,z,\cdot)\|_{\mcr{L}^{2}(U_1)}^2
  \leq  C|y|^2 + \alpha|z|^2;
 \eeqnn

 \item[(3)]
 for some constant $K>0$ we have
   $\beta(s,y,z,\zeta)=h(s,y,z)+\int_FC(s,u)\zeta_s(u)\nu(du)$ with $$|h(s,y,z)|\leq K(|y|+|z|),$$
   where $C(s,u)\geq -1$ and $\int_F|C(s,u)|^2\nu(du)\leq K$ for any $s\in [0,T]$.
  \end{enumerate}
   If $Y_T\leq 0$ a.s., we have $\mathbf{P}(Y_t \leq 0:t\in[0,T])=1.$
 \end{lemma}
 \proof
 Here we just prove $\mathbf{P}(Y_t \leq 0:t\in[0,T])=1$ under the corresponding conditions.
 It suffices to prove this theorem with $\nu(du)$ to be a finite Borel measure.
 Actually, for the general case we can always find a sequence $F_n\nearrow F$ such that $\nu_n(F)=\nu(F_n)<\infty$  and  $\nu_n(\cdot)=\mathbf{1}_{\{\cdot\in F_n\}}\nu(\cdot)\rightarrow \nu(\cdot).$
  For any $n\geq 1$, from Theorem~\ref{T1} and \ref{T2} there exists a unique solution $(Y^{(n)}_s,Z^{(n)}_s, \zeta^{(n)}_s(u))$ to (\ref{4.1}) with $M(ds,du)$ replaced by $M_n(ds,du)$, which has intensity $ds\nu_n(du)$.
  Like the proof of Theorem~\ref{T2}, we also have
  $Y^{(n)}_t\rightarrow Y_t\quad \mbox{in }\ \mathbb{S}^2_{\mathscr{F},T}.$
  For any integer $n\geq 0$, let
 $$\alpha_n=\exp\left\{-\frac{n(n+1)}{2}\right\}.$$
 Then $\alpha_n\rightarrow 0$ decreasingly as $n\rightarrow \infty$ and
 $$\int_{\alpha_n}^{\alpha_{n-1}}z^{-1}dz=n.$$
 Let $x\mapsto g_n(x)$ be a positive continuous function supported by $(\alpha_n,\alpha_{n-1})$ such that
 $$\int_{\alpha_n}^{\alpha_{n-1}}g_n(x)dx=1\qquad \mbox{and}\qquad xg_n(x)\leq\frac{2}{n}.$$
 Moreover, for any $n> 0$, define
 \beqnn
 f_n(z)=\Big|\int_0^{z}dy\int_0^yg_n(x)dx\Big|^2, \qquad z\in \mathbb{R}.
 \eeqnn
 It is easy to see that
 \begin{enumerate}
 \item[(a)] $f_n(z)\rightarrow |z^+|^2$ increasingly.

 \item[(b)] $
 f_n'(z)|=\left\{
 \begin{array}{ll}
 2\int_0^{z}g_n(x)dx\int_0^{z}dy\int_0^yg_n(x)dx\leq 2z, & z>0;\cr
 0, & z\leq 0
 \end{array}\right.
 $
 and $\lim\limits_{n\rightarrow\infty}f_n'(z)=2z^+$.
 \item[(c)] $
 f_n''(z)=\left\{
 \begin{array}{ll}
 2|\int_0^{z}g_n(x)dx|^2+2g_n(z)\int_0^{z}dy\int_0^yg_n(x)dx\leq 2+\frac{4}{n}, & z>0;\cr
 0, & z\leq 0
 \end{array}\right.
 $
 \vspace{10pt}\\
 and $\lim\limits_{n\rightarrow\infty}f_n''(z)=2\mathbf{1}_{\{z>0\}}$.
  \end{enumerate}
  Applying Proposition~\ref{P2.1} to $f_n(Y_t)$, Since $Y_T\leq 0$ a.s. we have
 \beqnn
 f_n(Y_t)\ar= \ar  \int_t^Tf_n'(Y_s)\beta(s,Y_s,Z_s,\zeta_s)ds+ \int_t^T\int_E f_n'(Y_s)\sigma(s,Y_s,Z_s,u)
 W^T(\overleftarrow{ds},du)\cr
 \ar\ar + \frac{1}{2} \int_t^Tds\int_Ef_n''(Y_s)|\sigma(s,Y_s,Z_s,u)|^2\pi(du)\cr
 \ar\ar
 +\int_{t-}^{T-}\int_{U_0} \big[f_n(Y_s+g_0(s,Y_s,Z_s,u)) - f_n(Y_s)\big]\tilde{N}_0^T(\overleftarrow{ds},du) \cr
 \ar\ar  +\int_t^Tds\int_{U_0} \big[f_n(Y_s+g_0(s,Y_s,Z_s,u)) - f_n(Y_s) - f_n'(Y_s)g_0(s,Y_s,Z_s,u) \big]\mu_0(du)\cr
 \ar\ar+\int_{t-}^{T-}\int_{U_1} \big[f_n(Y_s+g_1(s,Y_s,Z_s,u)) - f_n(Y_s)\big]N_1^T(\overleftarrow{ds},du) \cr
 \ar\ar
 - \int_t^Tf_n'(Y_s)Z_s dB_s- \frac{1}{2} \int_t^Tf_n''(Y_s)|Z_s|^2 ds \cr
 \ar\ar -\int_t^T\int_{F} \big[f_n(Y_s+\zeta_s(u)) - f_n(Y_s)\big]\tilde{M}(ds,du) \cr
 \ar\ar  -\int_t^Tds\int_{F} \big[f_n(Y_s+\zeta_s(u)) - f_n(Y_s)- f_n'(Y_s)\zeta_s(u) \big]\nu(du).
 \eeqnn
 Taking expectation to the above inequality, we have
 \beqnn
 \mathbf{E}\big[ f_n(Y_t)\big]\ar= \ar \int_t^T\mathbf{E}\big[f_n'(Y_s)\beta(s,Y_s,Z_s,\zeta_s)\big]ds +  \frac{1}{2}\int_t^Tds\int_E\mathbf{E}\big[f_n''(Y_s)|\sigma(s,Y_s,Z_s,u)|^2\big]\pi(du) \cr
  \ar\ar  +\int_t^Tds\int_{U_0} \mathbf{E}\big[f_n(Y_s+g_0(s,Y_s,Z_s,u)) - f_n(Y_s) - f_n'(Y_s)g_0(s,Y_s,Z_s,u) \big]\mu_0(du)\cr
  \ar\ar+\int_t^Tds\int_{U_1} \mathbf{E}\big[f_n(Y_s+g_1(s,Y_s,Z_s,u)) - f_n(Y_s) \big]\mu_1(du) -\frac{1}{2}\int_t^T\mathbf{E}\big[f_n''(Y_s)|Z_s|^2\big] ds\cr
  \ar\ar -\int_t^Tds\int_{F} \mathbf{E}\big[f_n(Y_s+\zeta_s(u)) - f_n(Y_s)- f_n'(Y_t)\zeta_s(u)\big]\nu(du).
 \eeqnn
 Since $Y_t\in\mathbb{S}^2_{\mathscr{F},T}$, from (a)-(c) and dominated convergence theorem, we have as $n\rightarrow \infty$
  \beqnn
 \mathbf{E}\big[|Y_t^+|^2\big]\ar\leq \ar \int_t^T\mathbf{E}\big[2Y_s^+\beta(s,Y_s,Z_s,\zeta_s)\big]ds +  \int_t^Tds\int_E\mathbf{E}\big[|\sigma(s,Y_s,Z_s,u)|^2\mathbf{1}_{\{Y_s>0\}}\big]\pi(du)\cr
  \ar\ar  +\int_t^Tds\int_{U_0} \mathbf{E}\big[|(Y_s+g_0(s,Y_s,Z_s,u))^+|^2 - |Y_s^+|^2 - 2Y_s^+g_0(s,Y_s,Z_s,u) \big]\mu_0(du)\cr
  \ar\ar+\int_t^Tds\int_{U_1} \mathbf{E}\big[|(Y_s+g_1(s,Y_s,Z_s,u))^+|^2 - |Y_s^+|^2\big]\mu_1(du)
  -\int_t^T\mathbf{E}\big[|Z_s|^2\mathbf{1}_{\{Y_s>0\}}\big] ds \cr
  \ar\ar  -\int_t^Tds\int_{F} \mathbf{E}\big[|(Y_s+\zeta_s(u))^+|^2 - |Y_s^+|^2- 2Y_s^+\zeta_s(u) \big]\nu(du).
 \eeqnn
 From condition (1) we have
 \beqlb
  \mathbf{E}\big[|Y_t^+|^2\big]\ar\leq \ar \int_t^T\mathbf{E}\big[2Y_s^+\beta(s,Y_s,Z_s,\zeta_s)\big] ds+  \int_t^Tds\int_E\mathbf{E}\big[|\sigma(s,Y_s,Z_s,u)|^2\mathbf{1}_{\{Y_s>0\}}\big]\pi(du)\cr
  \ar\ar+\int_t^Tds\int_{U_1} \mathbf{E}\big[2Y_s^+g_1(s,Y_s,Z_s,u)) + |g_1(s,Y_s,Z_s,u)|^2\mathbf{1}_{\{Y_s>0\}}\big]\mu_1(du)\cr
    \ar\ar  +\int_t^Tds\int_{U_0} \mathbf{E}\big[|g_0(s,Y_s,Z_s,u)|^2\mathbf{1}_{\{Y_s>0\}}\big]\mu_0(du) -  \int_t^T\mathbf{E}\big[|Z_s|^2\mathbf{1}_{\{Y_s>0\}}\big] ds \cr
  \ar\ar -\mathbf{E}\Big[\int_t^Tds\int_{F} [|(Y_s+\zeta_s(u))^+|^2 - |Y_s^+|^2- 2Y_s^+\zeta_s(u) ]\nu(du)\Big].\label{4.2}
 \eeqlb
 Let $\eta$ denote the integrand in the last term of (\ref{4.2}). Then
 \beqlb\label{4.3}
 \eta=\left\{
 \begin{array}{ll}
 |(Y_s+\zeta_s(u))^+|^2\geq 0, & \quad \mbox{if } Y_s\leq 0;\cr
 |\zeta_s(u)|^2  & \quad   \mbox{if } Y_s> 0 \mbox{ and } \zeta_s(u)\geq -Y_s;\cr
 - |Y_s|^2 - 2Y_s\zeta_s(u)  &  \quad   \mbox{if } Y_s> 0 \mbox{ and } \zeta_s(u)<-Y_s.
 \end{array}
 \right.
 \eeqlb
 Otherwise, by Cauchy's inequality and H\"{o}lder's inequality, for any $b>0$,
  \beqlb\label{4.4}
  2Y_s^+\int_FC(s,u)\zeta_s(u)\nu(du)\ar=\ar 2Y_s^+\int_FC(s,u)\zeta_s(u)\mathbf{1}_{\{\zeta_s(u)\geq -Y_s\}}\nu(du)\cr
  \ar\ar+2Y_s^+\int_FC(s,u)\zeta_s(u)\mathbf{1}_{\{\zeta_s(u)< -Y_s\}}\nu(du)\cr
  \ar \leq \ar b\left|\int_FC(s,u)\zeta_s(u)\mathbf{1}_{\{Y_s>0,\zeta_s(u)\geq -Y_s\}}\nu(du)\right|^2\cr
  \ar\ar+\frac{1}{b}|Y_s^+|^2+2Y_s^+\int_FC(s,u)\zeta_s(u)\mathbf{1}_{\{\zeta_s(u)< -Y_s\}}\nu(du)\cr
  \ar \leq \ar b\int_F|C(s,u)|^2\nu(du)\int_F|\zeta_s(u)|^2\mathbf{1}_{\{Y_s>0,\zeta_s(u)\geq -Y_s\}}\nu(du)\cr
  \ar\ar +\frac{1}{b}|Y_s^+|^2+2Y_s^+\int_FC(s,u)\zeta_s(u)\mathbf{1}_{\{\zeta_s(u)< -Y_s\}}\nu(du). \ \
  \eeqlb
 From (\ref{4.3}), (\ref{4.4}) and conditions in this theorem, for $a,c>0$ we have
 \beqnn
 \mathbf{E}\big[|Y_t^+|^2\big] \ar\leq \ar  \big(\frac{1}{a}+\frac{1}{b}+\frac{\mu_1(U_1)}{c}+\nu(F)\big)\int_t^T\mathbf{E}\big[|Y_s^+|^2\big]ds
 +a\int_t^T\mathbf{E}\big[|h(s,Y_s,Z_s)|^2\mathbf{1}_{\{Y_s>0\}}\big]ds \cr
  \ar\ar + \mathbf{E}\Big[\int_t^T\left(b\int_F|C(s,u)|^2\nu(du)-1\right)ds\int_F|\zeta_s(u)|^2\mathbf{1}_{\{Y_s>0,\zeta_s(u)\geq -Y_s\}}\nu(du)\Big]\cr
  \ar\ar + \mathbf{E}\Big[\int_t^T2Y_s^+ds\int_F[C(s,u)+1]\zeta_s(u)\mathbf{1}_{\{\zeta_s(u)< -Y_s\}}\nu(du)\Big]\cr
  \ar\ar+  \int_t^Tds\int_E\mathbf{E}\Big[|\sigma(s,Y_s,Z_s,u)|^2\mathbf{1}_{\{Y_s>0\}}\Big] \pi(du) \cr
  \ar\ar +\int_t^Tds\int_{U_0} \mathbf{E}\Big[|g_0(s,Y_s,Z_s,u)|^2 \mathbf{1}_{\{Y_s>0\}} \Big]\mu_0(du)- \int_t^T\mathbf{E}\Big[|Z_s|^2 \mathbf{1}_{\{Y_s>0\}}\Big]ds\cr
 \ar\ar+\int_t^Tds\int_{U_1} \mathbf{E}\Big[(1+c) |g_1(s,Y_s,Z_s,u)|^2\mathbf{1}_{\{Y_s>0\}}\Big]\mu_1(du).
 \eeqnn
  Here we choose $b$ small enough such that
 $b\int_F|C(s,u)|^2\nu(du) \leq 1$.
 From condition (2) and (3) we have
 \beqnn
 \mathbf{E}\big[ |Y_t^+|^2 \big]\ar\leq \ar  D\int_t^T\mathbf{E}\big[|Y_s^+|^2\big]ds +\tilde{\alpha}\int_t^T\mathbf{E}\Big[|Z_s|^2 \mathbf{1}_{\{Y_s>0\}}\Big]ds,
 \eeqnn
 where $D>0$ is a constant and $ \tilde{\alpha}=aK+(1+c)\alpha-1$. Let $a$ and $c$ be small enough such that $\tilde{\alpha}<0$, we have
 $$ \mathbf{E}\big[ |Y_t^+|^2 \big]\leq   D\int_t^T\mathbf{E}\big[|Y_s^+|^2\big]ds .$$
 By Gronwall's inequality, we have $ \mathbf{E}[ |Y_t^+|^2 ]=0$, which means
 $\mathbf{P}\left(Y_t\leq 0\right)=1$ for any $t\in[0,T]$.
 Like the proof of Theorem~\ref{T1}, we will get the desired result.
 \qed

 \begin{proposition}
 The conclusion in Lemma~\ref{Comparison.L1} remains true if (3) is replaced by
  \begin{enumerate}
  \item[(3')]   for some constant $K>0$ have $$|\beta(s,y,z,\zeta)|\leq K(|y|+|z|)+\int_FC(s,u)|\zeta(u)|\nu(du),$$
     where $0\leq C(s,u)\leq 1$ and $\int_F|C(s,u)|^2\nu(du)\leq K$ a.s. for any $s\in [0,T]$.
  \end{enumerate}
 \end{proposition}


 Now let us derive the general comparison theorem from Lemma~\ref{Comparison.L1}. However, since the deficiency of information about $Z_s$, we only consider the following case.
 For $i=1,2$, suppose $(Y_t^{(i)}, Z_t^{(i)},\zeta_s^{(i)}(u))$ is a solution to
 \beqlb\label{4.5}
 Y_t^{(i)}\ar=\ar Y_T^{(i)}+\int_t^T\beta^{(i)}(s,Y_s^{(i)},Z_s^{(i)},\zeta_s^{(i)})ds +
 \int_t^T\int_E\sigma(s,Y_s^{(i)},Z_s^{(i)},u)W^T(\overleftarrow{ds},du)\cr
 \ar\ar
 + \int_{t-}^{T-}\!\!\int_{U_0} g_0(s,Y_{s}^{(i)},u)
 \tilde{N}_0^T(\overleftarrow{ds},du)+ \int_{t-}^{T-}\!\!\int_{U_1}
 g_1(s,Y_{s}^{(i)},u) N_1^T(\overleftarrow{ds},du)\cr
 \ar\ar
 - \int_t^TZ_s^{(i)} dB_s- \int_t^T\!\!\int_F\zeta_s^{(i)}(u)\tilde{M}(ds,du).
 \eeqlb

 \begin{theorem}\label{T4.1}
 Suppose
 \begin{enumerate}
  \item[(1)] $\beta^{(1)}(s,y,z,\zeta)\leq \beta^{(2)}(s,y,z,\zeta)$;

  \item[(2)] both $y+g_0(s,y,u)$ and $y+g_1(s,y,u)$ are nondecreasing with respect to $y$;

  \item[(3)] $\sigma(s,y,z,u)$, $g_0(s,y,u)$ and $g_1(s,y,u)$ satisfy (\ref{2.3});


  \item[(4)] there exists a constant $K>0$ such that one of following conditions satisfies:
  \begin{enumerate}
   \item[(a)] $\beta^{(1)}(s,y,z,\zeta)=h^{(1)}(s,y,z)+\int_FC^{(1)}(s,u)\zeta_s(u)\nu(du)$ with
              $$|h^{(1)}(s,y,z)-h^{(1)}(s,y',z')|\leq K(|y-y'|+|z-z'|), $$
              where $C^{(1)}(s,u)\geq -1$ and $\int_F|C^{(1)}(s,u)|^2\nu(du)\leq K$ for any $s\in [0,T]$.
   \item[(b)] $\beta^{(2)}(s,y,z,\zeta)=h^{(2)}(s,y,z)+\int_FC^{(2)}(s,u)\zeta_s(u)\nu(du)$ with
              $$|h^{(2)}(s,y,z)-h^{(2)}(s,y',z')|\leq K(|y-y'|+|z-z'|),$$
             where $C^{(2)}(s,u)\leq 1$ and $\int_F|C^{(2)}(s,u)|^2\nu(du)\leq K$  for any $s\in [0,T]$.
  \end{enumerate}
 \end{enumerate}
 If $Y_T^{(1)}\leq Y_T^{(2)}$ a.s., then $\mathbf{P}(Y_t^{(1)}\leq Y_t^{(2)}:t\in[0,T])=1.$
 \end{theorem}
 \proof
 Here we assume (a) in condition (4) of this theorem holds and $\nu(du)$ to be a finite Borel measure. Let $(\bar{Y}_t, \bar{Z}_t, \bar{\zeta}_t(u)) = (Y_t^{(1)}-Y_t^{(2)},Z_t^{(1)}-Z_t^{(2)},\zeta_t^{(1)}(u)-\zeta_t^{(2)}(u))$. From (\ref{4.1}) and condition (1) we get
 \beqlb\label{4.6}
 \bar{Y}_t
 \ar\leq\ar \bar{Y}_T+
 \int_t^T\bar{\beta}(s)ds + \int_t^T\int_E\bar{\sigma}(s,u)
 W^T(\overleftarrow{ds},du)\cr
 \ar\ar
 +\int_{t-}^{T-}\int_{U_0} \bar{g_0}(s,u) \tilde{N}_0^T(\overleftarrow{ds},du) +
 \int_{t-}^{T-}\int_{U_1} \bar{g_1}(s,u) N_1^T(\overleftarrow{ds},du)\cr
 \ar\ar
 - \int_t^T\bar{Z}_s dB_s- \int_t^T\int_F\bar{\zeta}_s(u)\tilde{M}(ds,du),
 \eeqlb
 where $ \bar{Y}_T = Y^{(1)}_T-Y^{(2)}_T$ and
 \beqnn
 \bar{\beta}(s) \ar=\ar
 \beta^{(1)}(s,Y_s^{(1)},Z_s^{(1)},\zeta_s^{(1)})-\beta^{(1)}(s,Y_s^{(2)},Z_s^{(2)},\zeta_s^{(2)}),\cr
 \bar{\sigma}(s,u) \ar=\ar
 \sigma(s,Y_s^{(1)},Z_s^{(1)},u)-\sigma(s,Y_s^{(2)},Z_s^{(2)},u),\cr
 \bar{g_0}(s,u) \ar=\ar
 g_0(s,Y_s^{(1)},u)-g_0(s,Y_s^{(2)},u),\cr
 \bar{g_1}(s,u) \ar=\ar
 g_1(s,Y_s^{(1)},u)-g_1(s,Y_s^{(2)},u).
 \eeqnn
 It is easy to check that $\bar{\beta}(s)$, $\bar{\sigma}(s,u)$, $\bar{g}_0(s,u)$ and $\bar{g}_1(s,u)$ satisfy conditions in Lemma~\ref{Comparison.L1}, so the desired result follows.
 \qed

 \begin{remark}\label{R4.1}~
 \begin{enumerate}
 \item[(1)] Condition (2) means any jumps from $N_0$ and $N_1$ will not make $Y^{(1)}$ exceed
  $Y^{(2)}$.

 \item[(2)] Since we do not have enough information about $Z_s$ and $\zeta_s$, the result does not include the case of $g_0$ and $g_1$ depending on $Z$ and $\zeta$. This is still an open problem.

 \item[(3)] Obviously, Condition (4) can not be weakened to (\ref{2.2}) in Condition~\ref{C1}, counterexamples about BSDEs as a special case of BDSDEs can be seen in many works, such as Situ (2005, p.245).

 \item[(4)] Condition (4) can be replaced by
  \begin{enumerate}
  \item[(4')] one of following conditions satisfies
  $$|\beta^{(i)}(s,y,z,\zeta)-\beta^{(i)}(s,y',z',\zeta')|\leq K(|y-y'|+|z-z'|)+\int_FC(s,u)|\zeta(u)-\zeta'(u)|\nu(du),$$
  where $0\leq C(s,u)\leq 1$ and $\int_F|C(s,u)|^2\nu(du)\leq K$ for any $s\in [0,T]$.
  \end{enumerate}

 \end{enumerate}
 \end{remark}

 In the sequel of this section, we will show that comparison theorem still hold for the $1/2$-H\"older continuous case which is studied in He et al. (2014).  For $i=1,2$, suppose $(Y_t^{(i)},Z_t^{(i)},\zeta_t^{(i)})$ is a solution to the following BDSDE:
 \beqlb
 Y_t^{(i)}\ar=\ar Y_T^{(i)}+\int_t^T\beta^{(i)}(s,Y_s^{(i)},\zeta_s^{(i)})ds +
 \int_t^T\int_E\sigma(s,Y_s^{(i)},Z_s^{(i)},u)W^T(\overleftarrow{ds},du)\cr
 \ar\ar
 + \int_{t-}^{T-}\!\!\int_{U_0} g_0(s,Y_{s}^{(i)},u)
 \tilde{N}_0^T(\overleftarrow{ds},du)+ \int_{t-}^{T-}\!\!\int_{U_1}
 g_1(s,Y_{s}^{(i)},u) N_1^T(\overleftarrow{ds},du)\cr
 \ar\ar
 -\int_t^TZ_s^{(i)} dB_s-\int_t^T\!\!\int_F\zeta_s^{(i)}(u)\tilde{M}(ds,du).
 \eeqlb
 \begin{theorem}\label{T4.3}
 Suppose
 \begin{enumerate}
  \item[(1)] $\beta^{(1)}(s,y,\zeta)\leq \beta^{(2)}(s,y,\zeta)$ a.s..

  \item[(2)] both $y+g_0(s,y,u)$ and $y+g_1(s,y,u)$ are nondecreasing with respect to $y$.

  \item[(3)] For any $s\in[0,T]$ and $(y,z), (y',z')\in \mathbb{R}^2$,
  $$
 \int_{U_0}|g_0(s,y,u)-g_0(s,y',u)|^2\mu_0(du)+\int_{U_1}|g_1(s,y,u)-g_1(s,y',u)|\mu_1(du)\leq C|y-y'|
 $$
  and
  $$
  \int_E|\sigma(s,y,z,u)-\sigma(s,y',z',u)|^2\pi(du) \leq C|y-y'| +\alpha|z-z'|^2
  $$
where $C>0$ and $0\leq\alpha \leq 1$.

  \item[(4)]
  there exists a constant $K>0$ such that one of the following conditions is satisfied:
  \begin{enumerate}
   \item[(a)] $\beta^{(1)}(s,y,\zeta)=h^{(1)}(s,y)+\int_FC^{(1)}(s,u)\zeta_s(u)\nu(du)$ with
       $$|h^{(1)}(s,y)-h^{(1)}(s,y')|\leq K|y-y'|,$$
              where $C^{(1)}(s,u)\in[-1,0]$ and $\int_F|C^{(1)}(s,u)|\nu(du)\leq K$ for any $s\in[0,T]$.
   \item[(b)] $\beta^{(2)}(s,y,\zeta)=h^{(2)}(s,y)+\int_FC^{(2)}(s,u)\zeta_s(u)\nu(du)$ with
              $$|h^{(2)}(s,y)-h^{(2)}(s,y')|\leq K|y-y'|,$$
             where  $ C^{(2)}(s,u)\in[0,1]$ and $\int_F|C^{(2)}(s,u)|\nu(du)\leq K$ for any $s\in[0,T]$.
  \end{enumerate}
 \end{enumerate}
 If  $Y_T^{(1)}\leq Y_T^{(2)}$ a.s., we have
  $\mathbf{P}(Y_t^{(1)}\leq Y_t^{(2)}:t\in[0,T])=1.$

\end{theorem}
 \proof
 Here we just prove this theorem under condition $(a)$ holds. Moreover like the proof of Theorem~\ref{T4.1}, we assume $\nu(F)<\infty$. Let $(\bar{Y}_t, \bar{Z}_t,\bar{\zeta}_t(u)) = (Y_s^{(1)}-Y_s^{(2)},Z_s^{(1)}-Z_s^{(2)},\zeta_s^{(1)}-\zeta_s^{(2)})$.
  For any $n\geq 0$, recall $g_n(z)$ defined in the proof of Theorem~\ref{T4.1}, let
 \beqnn
 f_n(z)=\int_0^{z}dy\int_0^yg_n(x)dx, \qquad z\in \mathbb{R}.
 \eeqnn
 It is easy to see that
 \begin{enumerate}
 \item[(a)] $f_n(z)\rightarrow z^+$ increasingly.

 \item[(b)] $
 f_n'(z)=\left\{
 \begin{array}{ll}
 \int_0^{z^+}g_n(x)dx\leq 1, & z>0;\cr
 0, & z\leq 0
 \end{array}\right.
 $ and $\lim\limits_{n\rightarrow\infty}f_n'(z)\rightarrow \mathbf{1}_{\{z>0\}}$.
 \item[(c)] $
 zf_n''(z)=\left\{
 \begin{array}{ll}
 zg_n(z)\leq2/n, & z>0;\cr
 0, & z\leq 0.
 \end{array}\right.
 $

 \item[(d)] For any $az\geq 0$,
 \beqnn
  |f_n(a+z)- f_n(a)|=\Big|\int_{a^+}^{(a+z)^+}dy\int_0^yg_n(x)dx\Big|\leq |z|\mathbf{1}_{\{a>0\}}
 \eeqnn
 and
 \beqnn
 |f_n(a+z)- f_n(a)-zf_n'(a)|\leq \frac{1}{n}z^2.
 \eeqnn
  \end{enumerate}
 By Proposition~\ref{P2.1}, since $\bar{Y}_T\leq 0$ a.s., we have
 \beqnn
 f_n(\bar{Y}_t) \ar=\ar \int_t^T f_n'(\bar{Y}_s)\tilde{\beta}(s)ds+
 \int_t^T\int_E f_n'(\bar{Y}_s)\bar{\sigma}(s,u)W^T(\overleftarrow{ds},du)  \cr
 \ar\ar
  +\frac{1}{2}\int_t^T ds\int_Ef_n''(\bar{Y}_s)|\bar{\sigma}(s,u)|^2\pi(du) \cr
  \ar\ar
  +\int_{t-}^{T-}\int_{U_0}\big[f_n(\bar{Y}_{s}+\bar{g_0}(s,u))-f_n(\bar{Y}_{s})\big]\tilde{N}_0^T(\overleftarrow{ds},du)\cr
  \ar\ar+\int_t^Tds\int_{U_0}\big[f_n(\bar{Y}_s+\bar{g_0}(s,u))-f_n(\bar{Y}_s)-f_n'(\bar{Y}_s)\bar{g_0}(s,u)\big]\mu_0(du)\cr
 \ar\ar
  +\int_{t-}^{T-}\int_{U_1}\big[f_n(\bar{Y}_{s}+\bar{g_1}(s,u))-f_n(\bar{Y}_{s})\big]N_1^T(\overleftarrow{ds},du)\cr
 \ar\ar
   -\int_t^Tf_n'(\bar{Y}_s)\bar{Z}_s dB_s - \frac{1}{2}\int_t^Tf_n''(\bar{Y}_s)|\bar{Z}_s|^2 ds \cr
    \ar\ar
    -\int_t^T\int_F\big[f_n(\bar{Y}_{s}+\bar{\zeta}_s(u))-f_n(\bar{Y}_{s})\big]\tilde{M}(ds,du)\cr
    \ar\ar
-\int_t^T\int_F\big[f_n(\bar{Y}_s+\bar{\zeta}_s(u))-f_n(\bar{Y}_s)-f_n'(\bar{Y}_s)\bar{\zeta}_s(u)\big]\nu(du)ds,
 \eeqnn
 where $\tilde{\beta}(s)=\beta^{(1)}(s,Y_s^{(1)},\zeta_s^{(1)})-\beta^{(2)}(s,Y_s^{(2)},\zeta_s^{(2)})$ and $\bar{\sigma}(s,u)$, $\bar{g_0}(s,u)$, $\bar{g_1}(s,u)$ are defined like before.
 Since $f_n'(z),f_n''(z)\geq 0$ and $\beta^{(1)}(s,y,z)\leq \beta^{(2)}(s,y,z)$, we have
 \beqnn
 f_n(\bar{Y}_t) \ar\leq \ar \int_t^T f_n'(\bar{Y}_s)\bar{\beta}(s)ds+
 \int_t^T\int_E f_n'(\bar{Y}_s)\bar{\sigma}(s,u)W^T(\overleftarrow{ds},du)\cr
 \ar\ar +
  \frac{1}{2}\int_t^T f_n''(\bar{Y}_s)(C|\bar{Y}_s|+\alpha |\bar{Z}_s|^2)ds\cr
  \ar\ar
  +\int_t^T\int_{U_0}\big[f_n(\bar{Y}_{s}+\bar{g_0}(s,u))-f_n(\bar{Y}_{s})\big]\tilde{N}_0^T(\overleftarrow{ds},du)\cr
  \ar\ar+\int_t^Tds\int_{U_0}\big[f_n(\bar{Y}_s+\bar{g_0}(s,u))-f_n(\bar{Y}_s)-f_n'(\bar{Y}_s)\bar{g_0}(s,u)\big]\mu_0(du)\cr
 \ar\ar
  +\int_t^T\int_{U_1}\big[f_n(\bar{Y}_{s}+\bar{g_1}(s+,u))-f_n(\bar{Y}_{s})\big]N_1^T(\overleftarrow{ds},du) \cr
  \ar\ar-\int_t^Tf_n'(\bar{Y}_s)\bar{Z}_s dB_s  -\frac{1}{2}\int_t^Tf_n''(\bar{Y}_s)|\bar{Z}_s |^2ds\cr
 \ar\ar -\int_t^T\int_F\big[f_n(\bar{Y}_{s}+\bar{\zeta}_s(u))-f_n(\bar{Y}_{s})\big]\tilde{M}(ds,du)\cr
  \ar\ar -\int_t^Tds\int_F\big[f_n(\bar{Y}_s+\bar{\zeta}_s(u))-f_n(\bar{Y}_s)-f_n'(\bar{Y}_s)\bar{\zeta}_s(u)\big]\nu(du),
 \eeqnn
 where
 $$\bar{\beta}(s) =
 \beta^{(1)}(s,Y_s^{(1)},\zeta_s^{(1)})-\beta^{(1)}(s,Y_s^{(2)},\zeta_s^{(2)}).$$
 Taking the expectation, since $\alpha<1$, we have
  \beqnn
 \mathbf{E}\big[f_n(\bar{Y}_t)\big] \ar\leq \ar \int_t^T \mathbf{E}\big[f_n'(\bar{Y}_s)\bar{\beta}(s)\big]ds +
  \frac{C}{2}\int_t^T \mathbf{E}\Big[f_n''(\bar{Y}_s)|\bar{Y}_s|\Big]ds\cr
  \ar\ar+\int_t^Tds\int_{U_0}\mathbf{E}\Big[f_n(\bar{Y}_s+\bar{g_0}(s,u))-f_n(\bar{Y}_s)-f_n'(\bar{Y}_s)\bar{g_0}(s,u)\Big]\mu_0(du)\cr
 \ar\ar
  +\int_t^Tds\int_{U_1}\mathbf{E}\Big[f_n(\bar{Y}_s+\bar{g_1}(s,u))-f_n(\bar{Y}_s)\Big]\mu_1(du)\cr
  \ar\ar -\int_t^Tds\int_F\mathbf{E}\Big[f_n(\bar{Y}_s+\bar{\zeta}_s(u))-f_n(\bar{Y}_s)-f_n'(\bar{Y}_s)\bar{\zeta}_s(u)\Big]\nu(du).
 \eeqnn
 By the assumption of this theorem and $(c)$, we have
 $$ f_n''(\bar{Y}_s)|\bar{Y}_s|\rightarrow 0,\quad \mbox{as } n\rightarrow \infty.$$
 Moreover, from Lemma 3.1 in Li and Pu (2012),
 $$\int_{U_0}\big[f_n(\bar{Y}_s+\bar{g_0}(s,u))-f_n(\bar{Y}_s)-f_n'(\bar{Y}_s)\bar{g_0}(s,u)\big]\mu_0(du)\rightarrow 0,\quad \mbox{as } n\rightarrow \infty.$$
 Thus by Fatou's lemma, we have as $n\rightarrow \infty$,
 \beqlb\label{4.7}
 \mathbf{E}\big[\bar{Y}_t^+\big] \ar\leq \ar \int_t^T \mathbf{E}\left[\bar{\beta}_1(s)\mathbf{1}_{\{\bar{Y}_s>0\}}\right]ds
  +\int_t^Tds\int_{U_1}\mathbf{E}\big[(\bar{Y}_s+\bar{g_1}(s,u))^+ - \bar{Y}_s^+\big]\mu_1(du)\cr
  \ar\ar -\int_t^Tds\int_F\mathbf{E}\left[(\bar{Y}_s+\bar{\zeta}_s(u))^+ - \bar{Y}_s^+ - \bar{\zeta}_s(u)\mathbf{1}_{\{\bar{Y}_s>0\}}\right]\nu(du).
 \eeqlb
Now let us discuss the integrand of (\ref{4.7}) denoted by $\eta$:
 \beqlb\label{4.8}
 \eta=\left\{
 \begin{array}{ll}
 (\bar{Y}_s+\bar{\zeta}(s,u))^+\geq 0, & \quad \mbox{if } \bar{Y_s}\leq 0;\cr
 0,  & \quad   \mbox{if } \bar{Y_s}> 0 \mbox{ and } \bar{\zeta}(s,u))\geq -\bar{Y_s};\cr
 - \bar{Y}_s^+ -  \bar{\zeta}_s(u)\mathbf{1}_{\{\bar{Y}_s>0\}},  &  \quad   \mbox{if } \bar{Y_s}> 0 \mbox{ and } \bar{\zeta}(s,u))<-\bar{Y_s}.
 \end{array}
 \right.
 \eeqlb
 Moreover, since $-1\leq C(s,u)\leq 0$,
\beqnn
\int_FC(s,u)\bar{\zeta}_s(u)\nu(du)\ar=\ar \int_FC(s,u)\bar{\zeta}_s(u)\mathbf{1}_{\{\bar{\zeta}(s,u))\geq 0\}}\nu(du)\cr
\ar\ar+\int_FC(s,u)\bar{\zeta}_s(u)\mathbf{1}_{\{-\bar{Y_s}\leq\bar{\zeta}(s,u))< 0\}}\nu(du)\cr
\ar\ar+\int_FC(s,u)\bar{\zeta}_s(u)\mathbf{1}_{\{\bar{\zeta}(s,u))< -\bar{Y_s}< 0\}}\nu(du)\cr
\ar\leq\ar  |\bar{Y_s}|\int_F|C(s,u)|\nu(du)\cr
\ar\ar+\int_FC(s,u)\bar{\zeta}_s(u)\mathbf{1}_{\{\bar{\zeta}(s,u))< -\bar{Y_s}<0\}}\nu(du)\cr
\ar\leq\ar K|\bar{Y_s}|+\int_FC(s,u)\bar{\zeta}_s(u)\mathbf{1}_{\{\bar{\zeta}(s,u))< -\bar{Y_s}<0\}}\nu(du).
\eeqnn
From this, (\ref{4.7}) and conditions in this theorem we have
\beqnn
 \mathbf{E}\big[\bar{Y}_t^+\big] \ar\leq \ar [2K+C+\nu(F)]\int_t^T \mathbf{E}\big[\bar{Y_s}^+\big]ds+\int_F[C(s,u)+1]\bar{\zeta}_s(u)\mathbf{1}_{\{\bar{\zeta}(s,u))< -\bar{Y_s}<0\}}\nu(du)\cr
\ar\leq\ar  [2K+C+\nu(F)]\int_t^T \mathbf{E}\big[\bar{Y_s}^+\big]ds.
\eeqnn
 By Gronwall's inequality, we have $\mathbf{E}\big[\bar{Y}_t^+\big]=0$ and
 $\mathbf{P}\big(Y_t^{(1)}\leq Y_t^{(2)}\big)=1$ for any $t\in[0,T]$.
 Like the proof of Theorem~\ref{T1}, we will get the desired result.
 \qed

 \section{Maximum and Minimum Solutions}
 \setcounter{equation}{0}
 \medskip

  In the proof of Theorem~\ref{T2}, Picard iteration argument  seriously depends on the Lipschitz condition.
  Actually, sometimes solutions still exist (maybe not unique), even the drift term is linear increasing which is much weaker than Lipschitz condition.  As a simple application of comparison theorems, in this section we will prove the existence of solution to (\ref{4.5}) under some weak conditions.

\begin{theorem}\label{T5.1}
Suppose conditions~(2), (3) in Theorem~\ref{T4.1} holds and there exists a constant $K>0$ such that one of the following conditions holds:
\begin{enumerate}
   \item[(a)] $\beta(s,y,z,\zeta)=h(s,y,z)+\int_FC(s,u)\zeta_s(u)\nu(du)$ and
              $|h(s,y,z)|\leq K(1+|y|+|z|),$
              where $C(s,u)\in (-\infty, 1]$ (or $C(s,u)\in [-1,\infty)$) and $\int_F|C(s,u)|^2\nu(du)\leq K$ for any $s\in[0,T]$.
   \item[(b)] $|\beta(s,y,z,\zeta)|\leq K(1+|y|+|z|+\|\zeta\|_{\mcr{L}^2(F)})$ and for any $y,z$
$$|\beta(s,y,z,\zeta)-\beta(s,y,z,\zeta')|\leq \int_FC(s,u)|\zeta_s(u)-\zeta_s'(u)|\nu(du),$$
where $C(s,u)\in[0,1]$ and $\sup_{s\in[0,T]}\int_F|C(s,u)|^2\nu(du)\leq K$ a.s.
\end{enumerate}
Then solutions to (\ref{4.5}) exist. Moreover, there exist two solutions $(Y_t^{I},Z_t^{I},\zeta_t^{I})$ and $(Y_t^{S},Z_t^{S},\zeta_t^{S})$ such that for any solution $(Y_t,Z_t,\zeta_t)$ to (\ref{4.1}) have
$$\mathbf{P}\left(Y_t^{I}\leq Y_t\leq Y_t^{S}:t\in[0,T]\right)=1.$$
\end{theorem}
 Before using comparison theorem to prove this theorem, we need to construct a suitable sequence of BDSDEs with solutions exist and satisfy the conditions of comparison theorems; see the following lemma.
 Since the proof is easy and similar to Lemma~1 in Lepeltier and Martin (1997), we will omit it.
\begin{lemma}\label{L5.1}
For $n\geq K$, let
 $$\beta^{I}_n(s,y,z,\zeta)=\inf_{y',z'\in \mathbb{R}^2}\Big\{\beta(s,y',z',\zeta)+n|y-y'|+n|z-z'|\Big\}$$
and
$$\beta^{S}_n(s,y,z,\zeta)=\min\Big\{\beta(s,y,z,\zeta)+K,\sup_{y',z'\in \mathbb{R}^2}\left\{\beta(s,y',z',\zeta)-n|y-y'|-n|z-z'|\right\}\Big\}.$$
Then $\beta^{I}_n(s,y,z,\zeta)$ and $\beta^{S}_n(s,y,z,\zeta)$ are $\mathscr{F}_t^r$-progressive and satisfy:
\begin{enumerate}
\item[(1)]  For any $n \geq K $, $ \beta^{I}_n(s,y,z,\zeta)\leq \beta^{I}_{n+1}(s,y,z,\zeta)\leq \beta(s,y,z,\zeta)$ and
$$\beta(s,y,z,\zeta)\leq\beta^{S}_{(n+1)}(s,y,z,\zeta)\leq \beta^{S}_n(s,y,z,\zeta)\leq\beta(s,y,z,\zeta)+K.$$

\item[(2)] If $\beta(s,y,z,\zeta)$ satisfies (a)(or (b)) in Theorem~\ref{T5.1}, then so do $\beta^{I}_n(s,y,z,\zeta)$ and $\beta^{S}_n(s,y,z,\zeta)$.

\item[(3)] For any $(y,z),(y',z')\in\mathbb{R}^2$ have
\beqnn
|\beta^{I}_n(s,y,z,\zeta)-\beta^{I}_n(s,y',z',\zeta)|\ar\leq\ar n(|y-y'|+|z-z'|),\cr
\ar\ar\cr
|\beta^{S}_n(s,y,z,\zeta)-\beta^{S}_n(s,y',z',\zeta)|\ar\leq\ar n(|y-y'|+|z-z'|).
\eeqnn

\item[(4)] if $(y_n,z_n,\zeta_n)\rightarrow (y,z,\zeta)$, then $\beta^{I}_n(s,y_n,z_n,\zeta_n)$ and $\beta^{S}_n(s,y_n,z_n,\zeta_n)$ converge to $\beta(s,y,z,\zeta)$.

\end{enumerate}
\end{lemma}

  From Theorm~\ref{T1} and \ref{T2}, there exist unique solutions to (\ref{4.5}) with $\beta$ replaced by $\beta^{I}_n$  and $\beta^{S}_n$ respectively, denoted by $(Y^I_{n}(t),Z^I_{n}(t),\zeta^I_{n}(t))$ and $(Y^S_{n}(t),Z^S_{n}(t),\zeta^S_{n}(t))$ .
  According to this lemma and Theorem~\ref{T4.1}, for any $n\geq K$ we have $Y^I_{n}(t)\leq Y^I_{n+1}(t) \leq Y^S_{n+1}(t)\leq Y^S_{n}(t)$, which means both $\{Y^I_{n}(t)\}$ and $\{Y^S_{n}(t)\}$ are convergent in $\mcr{L}_{\mcr{F},T}^2$.
  So it suffices to show $(Y^I_{n}(t),Z^I_{n}(t),\zeta^I_{n}(t))\rightarrow (Y_t^{I},Z_t^{I},\zeta_t^{I})$ and $(Y^S_{n}(t),Z^S_{n}(t),\zeta^S_{n}(t))\rightarrow (Y_t^{S},Z_t^{S},\zeta_t^{S})$ in $\mcr{L}_{\mcr{F},T}^2\times \mcr{L}_{\mcr{F},T}^2\times \mcr{L}_{\mcr{F},T}^2(F)$ as $n\rightarrow \infty$. The key point is to prove for any $t\in[0,T]$,
 \beqnn
 \int_t^T\beta^{I}_n(s,Y^I_{n}(s),Z^I_{n}(s),\zeta^I_{n}(s))ds \ar\rightarrow\ar
 \int_t^T\beta(s,Y_s^{I},Z_s^{I},\zeta_s^{I})ds,\cr
 \ar\ar\cr
 \int_t^T\beta^{S}_n(s,Y^S_{n}(s),Z^S_{n}(s),\zeta^S_{n}(s))ds \ar\rightarrow\ar \int_t^T\beta(s,Y_s^{S},Z_s^{S},\zeta_s^{S})ds.
 \eeqnn
 By condition (4) in Lemma~\ref{L5.1}, it suffices to prove $(Y^I_{n}(t),Z^I_{n}(t),\zeta^I_{n}(t))$ and $(Y^S_{n}(t),Z^S_{n}(t),\zeta^S_{n}(t))$ are uniformly bounded on $[0,T]$.

\begin{lemma}\label{L5.2}
Assume conditions in Theorem~\ref{T5.1} hold. Then there exists $C>0$ such that for any $n\geq K$ such that
\beqnn
\|Y^I_{n}\|_{\mathbb{S}_{T}^2} \vee\|Y^S_{n}\|_{\mathbb{S}_{T}^2} \vee
\|Z^I_{n}\|_{\mcr{L}_T^2}\vee \|Z^S_{n}\|_{\mcr{L}_T^2}\vee\|\zeta^I_{n}\|_{\mcr{L}_T^2(E)}\vee \|\zeta^S_{n}\|_{\mcr{L}_T^2(E)}\leq C.
\eeqnn
\end{lemma}
 \proof Here we just prove this lemma with (a) in Theorem~\ref{T5.1} holds. First let us find a strip to cover all $Y^{I}_n(t)$ and $Y^{S}_n(t)$, which is easier to be dealt with. Define
 \beqnn
 \beta^{*}(s,y,z,\zeta)\ar=\ar K(2+|y|+|z|)+\int_FC(s,u)\zeta_s(u)\nu(du),\cr
 \beta_{*}(s,y,z,\zeta)\ar=\ar -K(2+|y|+|z|)+\int_FC(s,u)\zeta_s(u)\nu(du).
 \eeqnn
 Obviously, $\beta^{*}$ and $\beta_{*}$ satisfy (\ref{2.2}) in Condition~\ref{C1}. From Theorem~\ref{T1} and \ref{T2}, solutions to (\ref{4.1}) uniquely exist with $\beta$ replaced by $\beta^{*}$ and $\beta_{*}$ respectively, denoted by $(Y^{*}_t,Z^{*}_t,\zeta^{*}_t)$ and $(Y_{*t},Z_{*t},\zeta_{*t})$. Moreover, we have $\beta_{*}(t,y,z,\zeta)\leq\beta^{I}_n(t,y,z,\zeta) \leq \beta^{S}_n(t,y,z,\zeta)\leq \beta^{*}(t,y,z,\zeta)$ for any $n\geq K$. So $Y_{*t}\leq Y^{I}_n(t)\leq Y^{S}_n(t)\leq Y^{*}_t $.
 It suffices to prove
 $\|Y^{*}\|_{\mathbb{S}_{T}^2} \vee\|Y_{*}\|_{\mathbb{S}_{T}^2} \leq D$.
 By Proposition~\ref{P2.1},
 \beqnn
 \lefteqn{|Y^{*}_t|^2+\int_t^T|Z^{*}_s|^2ds+\int_t^T \|\zeta^{*}_s\|_{\mcr{L}^{2}(F)}^{2}ds = }\ar\ar\cr
 \ar\ar|Y_T|^2+ 2K\int_t^TY^{*}_s(2+ |Y^{*}_s|+|Z^{*}_s|)ds+2\int_t^TY^{*}_s ds\int_FC(s,u)\zeta^{*}_s(u)\nu(du)\cr
 \ar\ar+2\int_t^T\int_EY^{*}_s\sigma(s, Y^{*}_s,Z^{*}_s,u)W^T(\overleftarrow{ds},du)++ \int_t^T\|\sigma(s, Y^{*}_s,Z^{*}_s,\cdot)\|_{\mcr{L}^2(E)}^2ds\cr
 \ar\ar +\int_{t-}^{T-}\int_{U_0}\big[|Y^{*}_{s}+g_0(s, Y^{*}_{s},u)|^2-|Y^{*}_{s}|^2\big]\tilde{N}_0^T(\overleftarrow{ds},du)+\int_t^T\|g_0(s, Y^{*}_s,\cdot)\|_{\mcr{L}^2(U_0)}^2 ds\cr \ar\ar+\int_{t-}^{T-}\int_{U_1}\big[|Y^{*}_{s}+g_1(s, Y^{*}_{s},u)|^2-|Y^{*}_{s}|^2\big]N_1^T(\overleftarrow{ds},du)\cr
 \ar\ar-2\int_t^TY^{*}_sZ^{*}_sdB(s) -\int_t^T\int_{F}\big[|Y^{*}_{s}+\zeta^{*}_s(u)|^2-|Y^{*}_{s}|^2\big]\tilde{M}(ds,du)\cr
 \ar\leq \ar  |Y_T|^2+2T+ (3K+K/a+1/b+\mu_1(U_1))\int_t^T|Y^{*}_s|^2ds+a\int_t^T|Z^{*}_s|^2ds\cr
 \ar\ar+bK\int_t^T \|\zeta^{*}_s\|_{\mcr{L}^{2}(F)}^{2}ds+2\int_t^T\int_EY^{*}_s\sigma(s, Y^{*}_s,Z^{*}_s,u)W^T(\overleftarrow{ds},du)\cr
 \ar\ar+ (1+c)\int_t^T\big[C|Y^{*}_s|^2+\alpha|Z^{*}_s|^2\big]ds+(1+1/c)\int_t^T\|\sigma(s, 0,0,\cdot)\|_{\mcr{L}^{2}(E)}^{2}ds\cr
 \ar\ar+2\int_t^T\|g_0(s, Y^{*}_s,\cdot)-g_0(s, 0,\cdot)\|_{\mcr{L}^{2}(U_0)}^{2}ds+2\int_t^T\|g_0(s, 0,\cdot)\|_{\mcr{L}^{2}(U_0)}^{2}ds\cr
 \ar\ar+4\int_t^T\|g_1(s, Y^{*}_s,\cdot)-g_1(s, 0,\cdot)\|_{\mcr{L}^{2}(U_1)}^{2}ds+4\int_t^T\|g_1(s, 0,\cdot)|_{\mcr{L}^{2}(U_1)}^{2}ds\cr
 \ar\ar+\int_{t-}^{T-}\int_{U_0}\big[|Y^{*}_{s}+g_0(s, Y^{*}_{s},u)|^2-|Y^{*}_{s}|^2\big]\tilde{N}_0^T(\overleftarrow{ds},du) \cr
\ar\ar+\int_{t-}^{T-}\int_{U_1}\big[|Y^{*}_{s}+g_1(s, Y^{*}_{s},u)|^2-|Y^{*}_{s}|^2\big]\tilde{N}_1^T(\overleftarrow{ds},du)\cr
\ar\ar -2\int_t^TY^{*}_sZ^{*}_sdB(s)-\int_t^T\int_{F}\big[|Y^{*}_{s}+\zeta^{*}_s(u)|^2-|Y^{*}_{s}|^2\big]\tilde{M}(ds,du).
\eeqnn
  Apply Burkholder-Davis-Gundy Inequality to this formula, for example, there are some constants $A,B>0$ such that
\beqnn
\lefteqn{\mathbf{E}\Big[\sup_{t\in[0,T]}\int_t^T\int_{U_0}[|Y^{*}_s+g_0(s, Y^{*}_s,u)|^2-|Y^{*}_s|^2]\tilde{N}_0^T(\overleftarrow{ds},du)\Big]}\ar\ar\cr
\ar\leq\ar A \mathbf{E}\Big[2\Big[\int_0^T|Y^{*}_s|^2 \|g_0(s, Y^{*}_s,\cdot)\|_{\mcr{L}^{2}(U_0)}^{2} ds\Big]^{1/2}\Big]+ B\mathbf{E}\Big[\int_t^T\|g_0(s, Y^{*}_s,\cdot)\|_{\mcr{L}^{2}(U_0)}^{2}ds\Big]\cr
\ar\leq\ar A \mathbf{E}\Big[2\Big[\sup_{t\in[0,T]}|Y^{*}_s|^2\Big]^{1/2}\Big[\int_0^T\|g_0(s, Y^{*}_s,\cdot)\|_{\mcr{L}^{2}(U_0)}^{2}ds\Big]^{1/2}\Big]+ B\mathbf{E}\Big[\int_t^T\|g_0(s, Y^{*}_s,\cdot)\|_{\mcr{L}^{2}(U_0)}^{2}ds\Big]\cr
\ar\leq\ar d \|Y^{*}\|_{\mbb{S}_T^{2}}^2+ (A/d+B)\mathbf{E}\Big[\int_t^T\|g_0(s, Y^{*}_s,\cdot)\|_{\mcr{L}^{2}(U_0)}^{2}ds\Big].
\eeqnn
 The last inequality above comes from Cauchy's inequality. Like the proof before and choose $d$ small enough, we have
\beqnn
C_0\|Y^{*}\|_{\mbb{S}_T^{2}}^2+\|Z^{*}\|_{\mcr{L}_T^{2}}^{2}+\|\zeta^{*}\|_{\mcr{L}_T^{2}(F)}^{2} \leq  C_1+C_2\|Y^{*}\|_{\mcr{L}_T^{2}}^{2} +C_3 \big[\|Z^{*}\|_{\mcr{L}_T^{2}}^{2}+\|\zeta^{*}\|_{\mcr{L}_T^{2}(F)}^{2}\big],
\eeqnn
 where $C_1,C_2,C_3>0$ and $C_0\in(0,1)$. Since $(Y^{*}_t,Z^{*}_t,\zeta^{*}_t)\in\mcr{L}_{\mcr{F},T}^2\times \mcr{L}_{\mcr{F},T}^2\times \mcr{L}_{\mcr{F},T}^2(F)$, so there exists $C>0$ such that $\|Y^{*}\|_{\mbb{S}_T^{2}}^2\leq C$.
 Similarly,   $\|Y_{*}\|_{\mbb{S}_T^{2}}^2\leq C$ also can be proved.
 Here we have proved the first part of this lemma.
 For the second part, we just prove  $\|Z^I_{n}\|_{\mcr{L}_T^2}\vee \|\zeta^I_{n}\|_{\mcr{L}_T^2}\leq C$, others are similar. Apply  the It\^{o}-Pardoux-Peng formula to $|Y^{I}_n(t)|^2$ and like the proof of Theorem~\ref{T1}, we have
\beqnn
\mathbf{E}[|Y^{I}_{n}(0)|^2]+a\|Y^{I}_{n}\|_{\mcr{L}_T^2}^2+b\|\zeta_{n}^{I}\|^2_{\mcr{L}_T^2(F)}\leq C_0 + C_1 \|Y^{I}_n\|_{\mcr{L}_T^{2}}^2\leq C,
\eeqnn
where $a,b\in(0,1)$ and $C_0,C_1>0$ independent to $n$. Here we have finished the proof.
\qed

  Since we have showed $\{Y^{I}_n(t)\}$ and $\{Y^{S}_n(t)\}$ converge, it suffices to identify $(Z^{I}_n,\zeta^{I}_n)$ and $(Z^{S}_n,\zeta^{S}_n)$ are Cauchy sequences in $\mcr{L}_{\mcr{F},T}^2\times \mcr{L}_{\mcr{F},T}^2(F)$.
\begin{lemma}\label{L5.3}
Assume the conditions in Theorem~\ref{T5.1} holds, then both $(Z^{I}_n,\zeta^{I}_n)$ and $(Z^{S}_n,\zeta^{S}_n)$ are convergent in $\mcr{L}_{\mcr{F},T}^2\times \mcr{L}_{\mcr{F},T}^2(F)$.
\end{lemma}
\proof Like before we just prove  $(Z^{S}_n,\zeta^{S}_n)$ converges with condition (a) in Theorem~\ref{T5.1} holds. For any $n>m>K$, let $(Y^{S}_{n,m}(t), Z^{S}_{n,m}(t), \zeta^{S}_{n,m}(t,u))= (Y^{S}_n(t)-Y^{S}_{m}(t),Z^{S}_{n}(t)-Z^{S}_{m}(t), \zeta^{S}_{n}(t,u)-\zeta^{S}_{m}(t,u))$ which satisfies
 \beqnn
 Y^{S}_{n,m}(t)
 \ar=\ar
 \int_t^T\overline{\beta}^{(n,m)}(s)ds + \int_t^T\int_E\bar{\sigma}^{(n,m)}(s,u)
 W^T(\overleftarrow{ds},du)\cr
 \ar\ar
 +\int_t^T\int_{U_0} \bar{g_0}^{(n,m)}(s,u) \tilde{N}_0^T(\overleftarrow{ds},du) +
 \int_t^T\int_{U_1} \bar{g_1}^{(n,m)}(s,u) N_1^T(\overleftarrow{ds},du)\cr
 \ar\ar
 - \int_t^T Z_{n,m}^{S}(s) dB_s -
 \int_t^T\int_F\zeta^S_{n,m}(s,u)\tilde{M}(ds,du),
 \eeqnn
 where $\overline{\beta}^{(n,m)}(s)$, $\bar{\sigma}^{(n,m)}(s,u)$, $\bar{g_0}^{(n,m)}(s,u)$ and $\bar{g_1}^{(n,m)}(s,u)$ are defined like before.

 By Proposition~\ref{P2.1} and taking the expectation, we have
 \beqnn \|Z^{S}_{n,m}\|_{\mcr{L}_{T}^2}^2  + \|\zeta^{S}_{n,m}\|_{\mcr{L}_{T}^2}^2
 \ar\leq \ar
 \mathbf{E}\left[2\int_0^TY^{S}_{n,m}(s)\overline{\beta}^{(n,m)}(s)ds\right] +(3C+\mu_1(U_1))\|Y^{S}_{n,m}\|_{\mcr{L}_{T}^2}^2+\alpha\|Z^{S}_{n,m}\|_{\mcr{L}_{T}^2}^2.
 \eeqnn
  From Lemma~\ref{L5.2} we have $\sup_{s\in[0,T]}\mathbf{E}[|\bar{\beta}^{S}_{n,m}(s)|^2]<\infty$
  and
  \beqnn
  \mathbf{E}\left[2\int_0^TY^{S}_{n,m}(s)\overline{\beta}^{(n,m)}(s)ds\right]\ar\leq \ar 2\|Y^{S}_{n,m}\|_{\mcr{L}_{T}^2}^2\left(\mathbf{E}\Big[\int_0^T|\overline{\beta}^{(n,m)}(s)|^2ds\Big]\right)^{1/2}\leq 2\sqrt{CT}\|Y^{S}_{n,m}\|_{\mcr{L}_{\mcr{F},T}^2},
  \eeqnn
where $C>0$ is independent to $n$ and $m$. By H\"{o}lder inequality,
\beqnn
(1-\alpha) \|Z^{S}_{n,m}\|_{\mcr{L}_{\mcr{F},T}^2}^2+\|\zeta^{S}_{n,m}\|_{\mcr{L}_{\mcr{F},T}^2(F)}^2\leq  2\sqrt{CT}\|Y^{S}_{n,m}\|_{\mcr{L}_{\mcr{F},T}^2}+(3C+\mu_1(U_1))\|Y^{S}_{n,m}\|_{\mcr{L}_{\mcr{F},T}^2}^2
\eeqnn
Since $Y^{S}_{n}$ is a Cauchy sequence, so are $Z^{I}_n$ and $\zeta^{I}_n$.
\qed

\noindent\textit{Proof of Theorem~\ref{T5.1}.} With the preparations of Lemma~\ref{L5.1}, \ref{L5.2} and \ref{L5.3}, like the classical proof in SDE theory, we can get the desired result (omit the details).
\qed

\bigskip
\noindent\textbf{Acknowledgement:} I would like to thank Professor Zenghu Li for the enlightening discussion about (\ref{2.1}) and some comments on the proofs of comparison theorems.

 \noindent {\scshape Wei Xu}\\
 School of Mathematical Sciences,\\
  Beijing Normal University, \\
 Beijing 100875, People's Republic of China.\\
\textit{E-mails:}  xuwei@mail.bnu.edu.cn

 \end{document}